\documentclass[a4paper]{amsproc}
\usepackage{amssymb}
\usepackage{mathrsfs}
\usepackage{amsfonts}
\usepackage{amsthm}
\usepackage{amsmath,amscd}
\usepackage{amsxtra}     % Use various AMS packages
\usepackage{bm}
\usepackage[all,cmtip]{xy}
\usepackage{epsfig}
\usepackage{verbatim}
\usepackage{color}
\usepackage{soul}
\usepackage[shortlabels]{enumitem} %for allowing letters in enumerate
\usepackage{hyperref}
\usepackage{tikz}
\usetikzlibrary{arrows,matrix,shapes,trees}

\numberwithin{equation}{section}
\theoremstyle{plain}
 \newtheorem{thm}[equation]{Theorem} 
 \newtheorem*{theorem*}{Theorem}
 \newtheorem{prop}[equation]{Proposition}
 \newtheorem{lem}[equation]{Lemma} 
 \newtheorem{cor}[equation]{Corollary}

 \newtheorem{lem'}{``Lemma''}[section]

 \newtheorem{teorema}{Theorem}
  \newtheorem{corollario}[teorema]{Corollary}

\theoremstyle{definition}
 \newtheorem{ex}[equation]{Example} 
 \newtheorem{defn}[equation]{Definition} 
\theoremstyle{remark}
 \newtheorem{rmk}[equation]{Remark}

\newcommand{\mr}{\mathrm}
\newcommand{\mc}{\mathcal}
\newcommand{\ms}{\mathscr}

\newcommand{\N}{{\mathbb N}}
\newcommand{\Q}{{\mathbb Q}}
\newcommand{\Z}{{\mathbb Z}}

\newcommand{\F}{{\mathbb F}}
\newcommand{\Galk}{\mathcal{G}_{K}}
 
\newcommand{\Gm}{\mathbb{G}_{m}}
\newcommand{\Ga}{\mathbb{G}_a}
\newcommand{\Gml}{\mathbb{G}_{m,\mr{log}}}
 
\newcommand{\Gmlb}{({\mathbb{G}}_{\mr{m,log}}/\Gm)}

\newcommand{\Hone}{{H}^{\dagger}} % for $H_1$ in  BT groups with sst reduction
\newcommand{\Htwo}{{H}^{\ddagger}} % for $H_2$ in  BT groups with sst reduction
 
 \newcommand{\Nmono}{{\nu}}
 
\newcommand{\fs}{\mr{fs}}
\newcommand{\Sch}{\mr{Sch}}
\newcommand{\cl}{\mr{cl}}
\newcommand{\cris}{\mr{cris}}
\newcommand{\et}{\mr{\acute{e}t}}
\newcommand{\zar}{\mr{zar}}
\newcommand{\fl}{\mr{fl}}

\newcommand{\kfl}{\mr{kfl}}
\newcommand{\sh}{\mr{sh}}
\newcommand{\sst}{\mr{st}}
\DeclareMathOperator{\Spec}{Spec}

\newcommand{\Bcris}{B_{\mathrm{cris}}}

\newcommand{\bt}{\mathbf{BT}}
\newcommand{\cO}{{\ms O}}
\newcommand{\cM}{{\mathcal M}}
\newcommand{\rig}{\mathrm{rig}}

\newcommand{\cHom}{\mathcal{H}om} %sheaf Hom
\newcommand{\cExt}{\mathcal{E}xt} % sheaf Ext
\DeclareMathOperator{\Hom}{Hom}
\newcommand{\Ext}{\mathrm{Ext}}
\newcommand{\EXT}{\mathrm{EXT}}
\newcommand{\EXTc}{\mathrm{EXT}} %category of extensions
\newcommand{\HOMc}{\mathrm{HOM}} %category of hom
\newcommand{\Extpan}{\mathrm{Extpan}}
\newcommand{\EXTPAN}{\mathrm{EXTPAN}}
\newcommand{\gp}{\mathrm{gp}}
\newcommand{\Ab}{\mathbf{Ab}}
\newcommand{\fin}{\mathrm{fin}}
\DeclareMathOperator{\coker}{coker}

\newcommand\isomto{\stackrel{\sim}{\smash{\longrightarrow}\rule{0pt}{0.4ex}}}
 
\newcommand{\btlogr}{\mathbf{BT}_{S, {\mr r}}^{\log}}
\newcommand{\btlogd}{\mathbf{BT}_{S, {\mr d}}^{\log}}
\newcommand{\btlogc}{\mathbf{BT}_{S, {\mr c}}^{\log}}

\newcommand{\Rep}{\mathbf{Rep}}
\newcommand{\Repz}{\mathbf{Rep}_{{\mathbb Z}_p}}
\newcommand{\Repq}{\mathbf{Rep}_{{\mathbb Q}_p}}

\title[Log $p$-divisible groups and semistable representations]{Log $p$-divisible groups and semistable representations}



\author[Alessandra Bertapelle]{\bfseries Alessandra Bertapelle}
\author[Shanwen Wang]{\bfseries Shanwen Wang}
\author[Heer Zhao]{\bfseries Heer Zhao}

\address{
    Alessandra Bertapelle,
    Dipartimento di Matematica ``Tullio Levi-Civita'', 
    Universit\`a degli Studi di Padova, 
    Padova, Italy
    }
\email{alessandra.bertapelle@unipd.it}
\address{
    Shanwen Wang,
    School of mathematics, 
    Renmin University of China, 
    Beijing, 100872,
    China, 
     }
\email{s{\_}wang@ruc.edu.cn 
    }
    	
\address{
    Heer Zhao, 
    Institute for Advanced Study in Mathematics, 
    Harbin Institute of Technology, Harbin, 150001, 
    China, 
    	   }
\email{heer.zhao@gmail.com}

\begin{document}

\vspace{18mm} \setcounter{page}{1} \thispagestyle{empty}

\maketitle
\rightline{\`a la m\'emoire de Jean-Marc Fontaine}

\begin{abstract}
Let $\cO_K$ be a henselian DVR with field of fractions $K$ and residue field of characteristic $p>0$. Let $S$ denote $\Spec \cO_K$ endowed with the canonical log structure. We show that the generic fiber functor $\btlogd\to \bt^\sst_K$ between the category of dual representable log $p$-divisible groups over $S$ and the category of $p$-divisible groups with semistable reduction over $K$ is an equivalence. 
If $\cO_K$ is further complete with perfect residue field and of mixed characteristic, we show that $\btlogd$ is also equivalent to the category $\Repz^{\sst,\{0,1\}}(\Galk)$ of semistable Galois $\Z_p$-representations with Hodge-Tate weights in $\{0,1\}$. Finally, we show that the above equivalences are compatible with monodromies. Note that the monodromy in the category $\btlogd$ is integral, and thus compatibility provides an integral structure to the Fontaine monodromy associated to any object in $\Repz^{\sst,\{0,1\}}(\Galk)$.

\end{abstract}

\setcounter{tocdepth}{2}
\tableofcontents

 %%%%%%%%%%%%%%%%%%%%%%%%%
\section{Introduction}\label{Sec.intro}
%%%%%%%%%%%%%%%%%%%%%%%%%

 Let $\cO_K$ be a henselian discrete valuation ring with field of fractions $K$ and residue field $k$ of characteristic $p>0$. Let $W(k)$ be the ring of Witt vectors with coefficients in $k$ and $K_0$ its field of fractions. 
 %%even if k is non perfect, W(k) is an integral domain and a local ring.
For some results, we will need the following additional assumptions \begin{itemize}[($\ast$)]
    \item   $k$ is perfect and $K$ is a finite extension of $K_0$; in particular $\cO_K$ is of mixed characteristic, complete and totally ramified over $W(k)$.
\end{itemize} 

 One of the main scope of this paper is to extend Theorem \ref{classicalA} to the logarithmic context. For clarity, we recall some details.

\subsection{The classical case:  $p$-divisible groups and crystalline representations}\label{Subsection.pdiv-to-crysrep}
 Assume that $\cO_K$ satisfies $(\ast)$.
 Let $\overline K$ be a fixed algebraic closure of $K$ and $\Galk=\mathrm{Gal}(\overline K/K)$ the absolute Galois group of $K$. 
Let $\bt_{\cO_K}$ (resp. $\bt_K$) be the category of $p$-divisible groups (or Barsotti-Tate groups) over $\cO_K$ (resp. $K$). Let $\Repz(\Galk)$ (resp. $\Repq(\Galk)$) be the category of (continuous) representations of $\Galk$ on free $\mathbb{Z}_p$-modules of finite rank (resp. finite-dimensional $\Q_p$-vector spaces), and let   $\Rep_\bullet^{\cris,\{0,1\}}(\Galk)$ for $\bullet=\Z_p,\Q_p$,  denote the full subcategory of $\Rep_\bullet(\Galk)$  consisting of crystalline representations with Hodge-Tate weights in $\{0,1\}$. 

 For any $H_K=\varinjlim_n H_{K,n}$ in $\bt_K$, the \emph{Tate module} 
\[T_p(H_K)=\varprojlim_n H_{K,n}(\overline K)\]
of $H_K$ lies in $\Repz(\Galk)$ naturally, and $V_p(H_K):=T_p(H_K)\otimes_{\Z_p}\Q_p \in \Repq(\Galk)$. For $H\in\bt_{\cO_K}$, let 
$T_p(H):=T_p(H_K)$  and $V_p(H):=V_p(H_K)$. 
Then, one has the following well-known theorem.

\begin{teorema}[Breuil, Fontaine, Kisin, Raynaud, Tate] \label{classicalA}
The fully faithful functors
\[T_p\colon  \bt_{\cO_K}\to \Repz(\Galk), \quad    V_p\colon \bt_{\cO_K}\otimes\Q\to \Repq(\Galk) \]
 induce equivalences of categories
\[T_p\colon \bt_{\cO_K} \isomto \Repz^{\cris,\{0,1\}}(\Galk), \quad  V_p\colon \bt_{\cO_K}\otimes\Q\isomto \Repq^{\cris,\{0,1\}}(\Galk)  .\]
\end{teorema}
We refer to Section \ref{section.classical} for more details on the above result.

\subsection{The logarithmic case: log $p$-divisible groups and semistable representations}\label{Subsection.logpdiv-to-sstrep}
%%%%%%%%%%%%%%%%%%%%%%%%%
Let $\Rep_\bullet^{\sst}(\Galk)$ for $\bullet=\Z_p,\Q_p$,  denote the full subcategory of $\Rep_\bullet (\Galk)$  consisting of semistable representations.
In view of Theorem \ref{classicalA}, an object $T\in\Repz^{\cris,\{0,1\}}(\Galk)$ corresponds to a $p$-divisible group $H_K$ over $K$ having good reduction, i.e., $H_K$ extends to $H\in\bt_{\cO_K}$ which is unique up to unique isomorphism by Tate's theorem (see \cite[Thm. 4]{tat1}). It is natural to ask the following question
\begin{itemize}[$\bigstar$]
    \item Does an object $T\in\Repz^{\sst,\{0,1\}}(\Galk)$ correspond to a $p$-divisible group $H_K$ over $K$ having ``semistable reduction'' in whatever sense?
\end{itemize}

To answer the question $\bigstar$, one has to define the notion ``having semistable reduction'' for $p$-divisible groups over $K$ which is less clear than the notation ``having good reduction''. This notion is supplied by de Jong in \cite[Def. 2.2]{dJong98} (see Definition \ref{def.sstr}) without the assumption $(\ast)$.  We denote by $\bt_K^{\sst}$ the full subcategory of $\bt_K$ consisting of $p$-divisible groups having semistable reduction in the sense of de Jong.

However, de Jong's definition describes having semistable reduction without specifying the degeneration. This is simply because it is not possible to construct degenerations of $p$-divisible groups in the classical geometric world. It is well-known that log geometry is the perfect framework for dealing with degeneration. Let $S:=\Spec({\ms O}_K)$ endowed with the canonical log structure. Kato introduced log $p$-divisible groups in \cite{kat4} (see Definition \ref{def.logpdivr}), and dual representable log $p$-divisible groups (see \eqref{def.logpdivd}) over $S$ serve as the degeneration of $p$-divisible groups having semistable reduction over $K$. Let $\btlogd$ be the category of dual representable log $p$-divisible groups over $S$. 

In order to answer the question $\bigstar$ under the stronger assumption $(\ast)$, we study the relations among the three (instead of two in the crystalline case) categories $\Repz^{\sst,\{0,1\}}(\Galk)$, $\bt_K^{\sst}$, and $\btlogd$. For the relation between $\bt_K^{\sst}$ and $\btlogd$, we do not add the assumption $(\ast)$.

 Since the log structure of $S$ is supported on the closed point, the generic fiber $H_K:=H\times_S\Spec K$ of $H\in \btlogd$ is a classical $p$-divisible group over $K$ (see the proof of Theorem \ref{BZ} for more details). Thus, we have a natural functor
\[(\ )_K:\btlogd\to {\bf BT}_K, \quad H\mapsto H_K.\]
We denote the compositions 
\[\btlogd\xrightarrow{(\ )_K} {\bf BT}_{K}\xrightarrow{T_p}{\bf Rep}_{\Z_p}(\Galk)\]
and
\[\btlogd\otimes\Q\xrightarrow{(\ )_K} {\bf BT}_{K}\otimes\Q\xrightarrow{V_p}{\bf Rep}_{\Q_p}(\Galk)\] 
as $T_p$ and $V_p$,
respectively, by abuse of notation. Our first main result is the following theorem.

\begin{teorema}\label{main}
    \begin{enumerate}[(a)]
    \item The functor 
        \[(\ )_K\colon \btlogd\to {\bf BT}_{K}\]
        is fully faithful and has essential image the full subcategory $\bt_{K}^{\sst}$.
        \item Assume that $\cO_K$ satisfies the extra assumption $(\ast)$. Then the functor 
        \[T_p\colon \btlogd\to {\bf Rep}_{\Z_p}(\Galk)  \text{  (resp. $V_p\colon \btlogd\otimes\Q\to {\bf Rep}_{\Q_p}(\Galk)$)} \]
        is fully faithful and its essential image is exactly the full subcategory $\Repz^{\sst,\{0,1\}}(\Galk)$ (resp. $\Repq^{\sst,\{0,1\}}(\Galk)$).
    \end{enumerate}    
\end{teorema}
One can illustrate the above theorem through the following diagram of equivalences of categories
\begin{equation}\label{triangle_diagram_main3}
    \xymatrix{
&\btlogd  \ar[ld]_-{(\ )_K}^-{\simeq}\ar[rd]^-{T_p}_-{\simeq}  \\
{\bf BT}_{K}^{\mr{st}}\ar[rr]^-{\simeq}_-{T_p}  &&{\bf Rep}^{\mr{st},\{0,1\}}_{\Z_p}(\Galk).
}
\end{equation}
 Part (a) of the theorem follows from Theorem \ref{BZ} while part (b) is  Theorem \ref{BZb}.

 In Section \ref{section.corollaries} we recover from Theorem \ref{main}  the following $p$-adic N\'eron-Ogg-Shafarevich criterion for semistable reduction of abelian varieties.

 \begin{corollario}[Breuil, Kisin]\label{c.NOS}
 Assume that $\cO_K$ satisfies the extra assumption $(\ast)$, and let $A_K$ be an abelian variety over $K$. Then $A_K$ has semistable reduction if and only if $T_p(A_K)$ is a semistable Galois representation. 
 \end{corollario}

An anonymous referee brought our attention to the relation between log $p$-divisible groups and strongly divisible modules (defined in \cite[\S2.2]{Bre02}), as suggested by Breuil in  \cite[Rmk. 4.2.2.12]{Bre00}. We then noticed that another consequence of Theorem \ref{main}, also proved in Section \ref{section.corollaries}, is the following corollary.

\begin{corollario}\label{c.stronglydivisible}
 Assume that $\cO_K$ satisfies the extra assumption $(\ast)$, and $p\geq 3$. Then the category $\btlogd$ is anti-equivalent to the category of strongly divisible modules of weight $\leq 1$.
\end{corollario}
 
Strictly speaking, only part (b) of Theorem \ref{main} is the exact analogue of Theorem \ref{classicalA}.
One  could attempt to  adapt the proof strategy of Theorem \ref{classicalA} in order to establish Theorem \ref{main} (b). To this end, one needs to require at least
\begin{enumerate}
\item a logarithmic version of Tate's theorem \cite[Cor. 1]{tat1};
\item a logarithmic version of \cite[Cor. 2.2.6]{Kis06}, which is \cite[Thm. 1.4]{Bre00} when $p\neq2$ and $k\subset\overline{\F}_p$;
\item a logarithmic version of Raynaud's results \cite[\S2.1 and Prop. 2.3.1]{Ray74}.
\end{enumerate} 
To the authors' knowledge, these results have not appeared in the literature. The third author has an old unpublished draft on (1) in which a part of (3) has also been considered, but they are not enough to conclude.

Or, one could prove first Corollary \ref{c.stronglydivisible} and then deduce Theorem \ref{main} (b), at least for $p>2$. However, this seems a demanding job as well.

In our approach to Theorem \ref{main}, we use Fargues' results from \cite{Far19, Far22}, which do not require $p>2$, to relate the two categories $\bt_K^{\sst}$ and ${\bf Rep}^{\mr{st},\{0,1\}}_{\Z_p}(\Galk)$. Then we use Grothendieck's theory of panach\'ee extensions to study the relations between $\btlogd$ and $\bt_K^{\sst}$ with the help of Kummer log flat cohomology, and $K$ is only assumed to be henselian here (no completeness and no restriction on the residue field). The theory of panach\'ee extensions is also used to construct the Grothendieck monodromy for log $p$-divisible groups, see the next subsection. Log $p$-divisible groups are geometric objects defined in the Kummer log flat topology. We believe that our approach involving log geometry and Kummer log flat cohomology is very natural, although it might not be the only one. 

\subsection{Comparisons of monodromies}
%%%%%%%%%%%%%%%%%%%%%%%%%
Objects of the equivalent categories in \eqref{triangle_diagram_main3} are endowed with natural monodromy maps.  We complete our study of the  equivalences in \eqref{triangle_diagram_main3} by showing that they preserve monodromy.

Let $H\in\btlogd$, and let $0\to H^\circ\to H\to H^\et\to 0$ be the connected-\'etale decomposition of $H$ (see \cite[\S3.2]{w-z1}). By \cite[Prop. 3.9]{w-z1}, both $H^\circ$ and $H^\et$ are classical $p$-divisible groups. By Kato's classification theorem of log $p$-divisible group, see Theorem \ref{thm1.2}, any object $H\in\btlogd$ corresponds to a pair $(H^\cl,\beta)$, where 
\begin{itemize}
    \item $H^\cl$ is a classical $p$-divisible group and called the classical part,
    \item $\beta\colon H^\et(1)\to H^\circ$ is a homomorphism of classical $p$-divisible groups called the \emph{Kato monodromy map} of $H$.
\end{itemize}
Let $H^\mu$ be the multiplicative part of $H^\circ$, then $\beta$ factors through $H^\mu$ and we denote the resulting map $H^\et(1)\to H^\mu$ still as $\beta$, by abuse of notation. The map $\beta\colon H^\et(1)\to H^\mu$ corresponds to a pairing 
\[c(H)\colon T_p(H^\et)\otimes T_p((H^\mu)^*)\to\Z_p,\] which we call the \emph{Kato monodromy pairing} of $H$ (see Definition \ref{p-adic Kato pairing}).

Now consider the generic fiber $H_K$ of $H$. Grothendieck's theory of panach\'ee extensions furnishes a pairing 
\[c^{\mathrm{Gr}}(H_K)\colon T_p(H^\et)\otimes T_p((H^\mu)^*)\to\Z_p\]
to $H_K$ or equivalently to $H$, which we call the \emph{Grothendieck monodromy pairing} of $H_K$ or equivalently of $H$ (see Definition \ref{p-adic Kato pairing}).

Our next main theorem compares $c(H)$ with $c^{\mathrm{Gr}}(H)$ (see Theorem \ref{thm.monodromypairings}).

\begin{teorema}\label{t.cGr}
    For any $H\in\btlogd$ we have $c(H)=c^{\mathrm{Gr}}(H_K)$.
\end{teorema}

Assume further that $\cO_K$ satisfies the assumption $(\ast)$. In the equivalence triangle \eqref{triangle_diagram_main3}, it seems that there is no monodromy associated to the objects of $\Repz^{\sst,\{0,1\}}(\Galk)$. However, consider the functors
\[\xymatrix{\Repz^{\sst,\{0,1\}}(\Galk)\ar[r]_-{V(\ )} &
\Repq^{\sst,\{0,1\}}(\Galk)\ar[r]^-{\simeq}_-{D_{\sst}} &\underline{M}^{\mathrm{a},\{-1,0\}}
},\]
where 
\begin{itemize}
    \item the functor $V(\ )$ is given by $T\mapsto T\otimes_{\Z_p}\Q_p $,
    \item $\underline{M}^{\mathrm{a}}$ denotes the category of admissible filtered $(\varphi,N)$-modules over $K$ (see \cite[\S 4.1]{CF00}), and $\underline{M}^{\mathrm{a},\{-1,0\}}$ denotes the full subcategory of $\underline{M}^{\mathrm{a}}$ consisting of objects $D$ such that $\mathrm{Fil}^{-1}D_K=D_K$ and $\mathrm{Fil}^{1}D_K=0$,
    \item the functor $D_{\sst}\colon \Repq^{\sst}(\Galk)\to \underline{M}^{\mathrm{a}}$ associating $D_{\sst}(V)$ from \eqref{functor.Dst} to $V\in\Repq(\Galk)$ is an exact tensor equivalence by \cite[Prop. 4.2]{CF00} and it clearly restricts to an equivalence $\Repq^{\sst,\{0,1\}}(\Galk)\xrightarrow{\simeq} \underline{M}^{\mathrm{a},\{-1,0\}}$.
\end{itemize}
For any $T\in \Repz^{\sst,\{0,1\}}(\Galk)$, if we pass to $\underline{M}^{\mathrm{a},\{-1,0\}}$ along the two functors $V(\ )$ and $D_{\sst}$, we have the $K_0$-linear endomorphism $N$ on the $K_0$-vector space $D_{\sst}(V(T))$ and we call $N$ the \emph{Fontaine monodromy map} of $T$. 

Given $H\in\btlogd$, our last main theorem compares the Kato monodromy map $\beta$ of $H$ with the Fontaine monodromy map of $V_p(H)=V(T_p(H))$ (see Theorem \ref{thm.monodromyKatoFontaine} and the paragraph after its proof).
\begin{teorema}\label{t.KFmonodromy}
    Assume that $\cO_K$ satisfies the assumption $(\ast)$ and let $H\in\btlogd$. Then the Kato monodromy map $\beta$ of $H$ determines the Fontaine monodromy map of $T_p(H)$, and vice versa rationally.
\end{teorema}

The paper is structured as follows. 
In Section \ref{section.classical} we rapidly introduce notations and cite results in the literature that lead to the proof of Theorem \ref{classicalA}. Before generalizing this result to the logarithmic context, we recall in Section \ref{s.general_results} results on log $p$-divisible groups. A useful technical lemma on panach\'ee extensions is proved in Section \ref{s.obstruction}. The proof of Theorem \ref{main} starts with Section \ref{s.lpdsst} and ends in Section \ref{s.Bb}.
Corollaries \ref{c.NOS} and \ref{c.stronglydivisible} are proved in Section \ref{section.corollaries}. Finally Theorems \ref{t.cGr} and \ref{t.KFmonodromy} on monodromy are proved in Sections \ref{s.monodromy_pairings} and \ref{s.compatibilityKF}, respectively.

%%%%%%%%%%%%%%%%%%%%%%%%%
\section{Details on the classical case}\label{section.classical}
%%%%%%%%%%%%%%%%%%%%%%%%%

In this section we assume that $\cO_K$ satisfies $(\ast)$, and give some details on Theorem \ref{classicalA} and its proof. Also, we introduce notation used later. 

For $H\in\bt_{\cO_K}$, let $T_p(H):=T_p(H_K)$ and $V_p(H):=V_p(H_K)$.
Since $K$ is of characteristic 0, the functors
\[T_p\colon \bt_K\to \Repz(\Galk) \ \text{ and } \ V_p\colon \bt_K\otimes\Q\to \Repq(\Galk)\]
are equivalences of categories \cite[p. 725]{Far19}. Then by Tate's Theorem \cite[Thm. 4]{tat1} the functors 
\[T_p\colon \bt_{\cO_K}\to \Repz(\Galk) \text{ and } V_p\colon \bt_{\cO_K}\otimes\Q\to \Repq(\Galk)\]
are fully faithful. It is a natural question to ask what the essential images of the last two functors are.
The answer needs Fontaine's period rings. 

Let $\Bcris$ be the ring of crystalline periods and let $B_{\rm st}=\Bcris[u]$ be the ring of log-crystalline periods, where $u=\log[p^{\flat}]$ with $p^{\flat}=(p, p^{1/p}, \cdots)$ (see \cite{CF00}). Both rings are endowed with an action of Frobenius $\varphi$ and a decreasing filtration. Furthermore, $B_{\rm st}$ is endowed with a unique $\Bcris$-derivation $N$ such that $N(u)=-1$, called the \emph{monodromy operator}. For any $T\in {\bf Rep}_{\mathbb{Z}_p}(\Galk)$, set 
\begin{equation}\label{functor.QpRep}
    V(T):=T\otimes_{\Z_p}\Q_p,
\end{equation}
and consider the filtered $\varphi$-module (resp. filtered $(\varphi, N)$-module) 
\begin{equation}\label{functor.Dcris}
    D_{\cris}(T):=D_{\cris}(V(T)):=(V(T)\otimes_{\Q_p} \Bcris)^{\Galk},
\end{equation}
\begin{equation}\label{functor.Dst}
    (\text{resp. } D_{\sst}(T):=D_{\sst}(V(T)):=(V(T)\otimes_{\Q_p} B_{\sst})^{\Galk} ).
\end{equation}
A $\mathbb{Z}_p$-representation $T$ of $\Galk$ is called \emph{crystalline} (resp. \emph{semistable}) if 
\[\dim_{\mathbb{Q}_p}V(T)=\dim_{K_0}D_{\cris}(T) \quad (\text{ resp. } \dim_{\mathbb{Q}_p}V(T)=\dim_{K_0} D_{\rm st}(T)),
\]  
and we denote the full subcategory of $\Repz(\Galk)$ consisting of crystalline (resp. semistable) representations by $\Repz^{\cris}(\Galk)$ (resp. $\Repz^{\sst}(\Galk)$). Note that a crystalline representation is automatically semistable, in other words $\Repz^{\cris}(\Galk)$ is a full subcategory of $\Repz^{\sst}(\Galk)$. 
Let $\bullet=\Z_p, \Q_p$. We have the following diagram of subcategories
\[\renewcommand{\arraystretch}{1.3}
\begin{array}[c]{ccccc}
\Rep_\bullet^{\sst,\{0,1\}}(\Galk) &\subset &\Rep_\bullet^{\sst}(\Galk) &\subset &\Rep_\bullet(\Galk)  \\
\rotatebox{90}{$\subset$}&&\rotatebox{90}{$\subset$}&\\
\Rep_\bullet^{\cris,\{0,1\}}(\Galk) &\subset &\Rep_\bullet^{\cris}(\Galk) 
\end{array},\]
where $\Rep_\bullet^{\cris,\{0,1\}}(\Galk)$ (resp. $\Rep_\bullet^{\sst,\{0,1\}}(\Galk)$) denotes the full subcategory of $\Rep_\bullet^{\cris}(\Galk)$ (resp. $\Rep_\bullet^{\sst}(\Galk)$) consisting of objects with Hodge-Tate weights in $\{0,1\}$.

Theorem \ref{classicalA} in the Introduction section answers the aforementioned question. 
For $p=2$, the theorem is simply a reformulation of \cite[Thm. 2.2.1]{Liu}, attributed in loc. cit. to Fontaine, Kisin, Raynaud and Tate. For completeness, we present below a full proof of Theorem \ref{classicalA}. It imitates the proof of \cite[Thm. 2.2.1]{Liu} and basically collects some famous results from the literature.

\emph{Proof of Theorem \ref{classicalA}}.
We only need to treat the case of $T_p$. By Tate's theorem \cite[\S4.2, Cor. 1]{tat1}, the functor $T_p$ is fully faithful. By Fontaine's theorem \cite[\S6.2]{Fon82}, the image of the functor $T_p$ lands in $\Repz^{\cris}(\Galk)$. It further lands in $\Repz^{\cris,\{0,1\}}(\Galk)$ by Tate's theorem \cite[\S4, Cor. 2 on p. 180]{tat1}. We are left to show that the essential image of $T_p$ is $\Repz^{\cris,\{0,1\}}(\Galk)$. 

Let $T$ be an object in $\Repz^{\cris,\{0,1\}}(\Galk)$. By Raynaud's theorem \cite[Prop. 2.3.1]{Ray74}, it suffices to show that $T/p^nT$ regarded as a finite group scheme over $K$ extends to a finite flat group scheme over $\cO_K$ for each $n$. Note that there exists $H\in \bt_{\cO_K}$ such that $T\otimes_{\Z_p}\Q_p\cong V_p(H)$, and we regard $T$ as a lattice inside $V_p(H)$ with respect to this isomorphism. The existence of such $H$ is due to Breuil when $p\geq3$ and $k\subset\overline{\F}_p$ (see \cite[Thm. 1.4]{Bre00}), and due to Kisin for the rest cases (see \cite[Cor. 2.2.6]{Kis06}). We are going to use $H$ to produce the aforementioned integral models. Let $a,b\in\N$ such that $p^{a+b}T_p(H)\subset p^aT\subset T_p(H)$. Then we have a diagram of finite $\Galk$-modules, that is, finite \'etale group schemes over $K$,
\[T/p^nT\cong p^aT/p^{a+n}T\xrightarrow{\alpha} T_p(H)/p^{a+n}T\xleftarrow{\beta} T_p(H)/p^{a+b+n}T_p(H)\cong (H_{a+b+n})_K\]
with $\alpha$ (resp. $\beta$) injective (resp. surjective). By \cite[\S2.1]{Ray74}, $\ker(\beta)$ extends to a finite flat group scheme $N$ over $\cO_K$, and thus the finite flat group scheme $H_{a+b+n}/N$ extends $T_p(H)/p^{a+n}T$. Applying \cite[\S2.1]{Ray74} once again, we see that $T/p^nT$ also extends to a finite flat group scheme over $\cO_K$. We are done.

%%%%%%%%%%%%%%%%%%%%%%%%%
\section{General results}\label{s.general_results}
%%%%%%%%%%%%%%%%%%%%%%%%%
This section mainly contains known results that will be used in the next sections. 

\subsection{Log $p$-divisible groups}\label{divcon}
In this subsection, we introduce Kato's theory of log $p$-divisible groups, which is developed in \cite{kat4}.
In the following, log structures are defined by sheaves of monoids for the \'etale topology.

Let $S$ be an fs log scheme whose underlying scheme $\mathring{S}$ is locally noetherian, and let $(\fs/S)$ be the category of fs log schemes over $S$. We endow $(\fs/S)$ with the Kummer log flat topology (cf. \cite[\S 2]{Kato21} and \cite[\S 2]{Niz08}), and denote the resulting site by $(\fs/S)_{\kfl}$.
Sometimes, we abbreviate $(\fs/S)_{\kfl}$ as $S_{\kfl}$ to shorten the formulas.  Similarly,  $(\fs/S)_{\fl}$ denotes the category $(\fs/S)$ with the classical flat topology (fppf). 
These two sites are denoted as $S_\fl^{\log}$ and $S_\fl^{\cl}$ respectively in \cite{kat4}.

\begin{defn}
Let $\Ab_{\kfl}(S)$ denote the category of sheaves of abelian groups over $(\mr{fs}/S)_{\kfl}$.
We define $(\mr{fin}/S)_{\mr{r}}$ as the full subcategory of $\Ab_{\kfl}(S)$ consisting of objects $F$ which are representable by an fs log scheme $f\colon F\to S$ such that the structure morphism $f$ is Kummer log flat and the underlying map of schemes is finite.
We call an object of $(\mr{fin}/S)_{\mr{r}}$ a \emph{finite Kummer log flat group log scheme}, or simply a \emph{finite kfl group log scheme}.
\end{defn}
Note that $F\in (\mr{fin}/S)_{\mr{r}}$ with $F\to S$ strict is just a classical finite flat group scheme over $\mathring{S}$ endowed with the log structure induced from $S$. We denote the full subcategory consisting of such objects by $(\mr{fin}/S)_{\mr{c}}$. Let $\Gm$ be the multiplicative group endowed with the induced log structure. For $F\in (\mr{fin}/S)_{\mr{r}}$, the \emph{Cartier dual} of $F$ is the sheaf \[F^*:=\mc{H}om_{S_{\mr{kfl}}}(F,\Gm).\]
The category $(\mr{fin}/S)_{\mr{d}}$ is the full subcategory of $(\mr{fin}/S)_{\mr{r}}$ consisting of objects $F$ with $F^*\in(\mr{fin}/S)_{\mr{r}}$.
 
\begin{defn}\label{def.logpdivr} 
A \emph{log $p$-divisible group} (or a \emph{log Barsotti-Tate group}) over $S$ is an object $H$ of $\Ab_{\kfl}(S)$ satisfying:
\begin{enumerate}[(a)]
\item $H=\varinjlim_nH_n$ with
$H_n:=\mr{ker}(p^n\colon H\rightarrow H)$; 
\item $p\colon H\rightarrow H$ is surjective;
\item $H_n\in (\mr{fin}/S)_{\mr{r}}$ for any $n> 0$.
\end{enumerate}
\end{defn}
We denote the category of log $p$-divisible groups over $S$ by $\btlogr$. We define full subcategories
\begin{equation}\label{def.logpdivd}
    \btlogc\subseteq \btlogd \subseteq \btlogr
\end{equation}
by: $H\in \btlogd $ (resp. $H\in \btlogc$), if $H_n\in (\mr{fin}/S)_{\mr{d}}$ (resp. $H_n\in (\mr{fin}/S)_{\mr{c}}$) for $n\geq1$.  We call the objects of $\btlogd$ the \emph{dual representable log $p$-divisible groups}. 
Clearly $H\in \btlogc$ amounts to a classical $p$-divisible group over $\mathring{S}$.

Recall that given a short exact sequence $0\to H'\to H\to H''\to 0$ in $\Ab_{\kfl}(S)$, if $H', H''$  are dual representable log $p$-divisible groups, the same is $H$ (cf. \cite[Prop. 2.3]{kat4}). 
 Furthermore, the exactness of a sequence of log $p$-divisible groups is equivalent to the exactness of the sequences of kernels of multiplication by $p^n$ for one (equivalently all) $n>0$.  

\subsection{Kato monodromy}\label{s.katomonodromy}
%%%%%%%%%%%%%%%%%%%%%%%%%
In this subsection, we assume that the underlying scheme of $S$ is $\mathring{S}=\Spec(A)$ with $A$ a noetherian henselian local ring. 
Let $k=A/\mathfrak{m}_A$ and $p:=\mr{char}(k)>0$. Suppose further that the log structure $\cM_S$ of $S$ admits a chart $\gamma\colon P_S\to \cM_S$ with $P$ an fs monoid, and that $\gamma$ induces an isomorphism $P=P_{S,\bar s} \xrightarrow{\sim}\cM_{S,\bar s}/\cO^\times_{S,\bar s} $,  where $\bar s$ denotes a geometric point above the closed point $s$ of $\Spec(A)$ and thus $\cO^\times_{S,\bar s} =(A^{\sh})^\times$, with $A^{\sh}$ the strict henselization of $A$.

\begin{rmk}\label{rem.stalks-charts}
 Note that the assumption on the chart $\gamma$ implies that the canonical map $\cM_{S,s}/\cO^\times_{S,s}\xrightarrow{\sim}\cM_{S,\bar s}/\cO^\times_{S,\bar s}$ is an isomorphism, where  $(-)_{s}$ denotes the stalk for the Zariski topology at $s$ and by abuse of notation $\cM_S$ also denotes the restriction $\cM_{S,\zar}$ of $\cM_S$ to the small Zariski site of $S$. 
 In fact, let $\eta\colon   S_{\et}\to   S_{\zar}$ denote the canonical map of small sites;
 the existence of the global chart $\gamma$ implies that $\eta^*\cM_{S,\zar} \xrightarrow{\sim}\cM_S$ by \cite[Ch. III, Prop. 1.4.1.2]{Ogu18}.
 Then for any \'etale neighborhood $(U,u)$ of $s$, the restriction of $\cM_S$ to the small Zariski site of $U$ is just the inverse image of $\cM_{S,\zar}$ along $U\to S$, therefore $\cM_{S,s}/\cO^\times_{S,s}\xrightarrow{\sim}\cM_{U,u}/\cO^\times_{U,u}$ by \cite[Ch. III, Rmk. 1.1.6]{Ogu18}. In particular, we get $\cM_{S,s}/\cO^\times_{S,s}\xrightarrow{\sim}\cM_{S,\bar s}/\cO^\times_{S,\bar s}$. Therefore, the requirement $P\xrightarrow{\sim}\cM_{S,\bar s}/\cO^\times_{S,\bar s}$ is equivalent to requiring $P \xrightarrow{\sim}\cM_{S,s}/\cO^\times_{S,s}$.
\end{rmk}

The torsion subgroups of an object of $\btlogd $ lie in $(\mr{fin}/S)_{\mr{d}}$, and the following theorem of Kato describes an object $F\in (\mr{fin}/S)_{\mr{d}}$ as an extension of classical finite flat group schemes.

\begin{prop}[Kato, \cite{kat4}] 
 Let $F\in (\mr{fin}/S)_{\mr{r}}$ and let $F^\circ$ be the connected component of $F$  that contains the image of the identity section. Then
\begin{enumerate}[(a)]
\item $F^\circ\in (\mr{fin}/S)_{\mr{c}}$.
\item $F^\et:=F/F^\circ\in (\mr{fin}/S)_{\mr{r}}$. 
\item Assume that $F$ is killed by a power of $p$. Then $F\in (\mr{fin}/S)_{\mr{d}}$ if and only if $F^\et\in (\mr{fin}/S)_{\mr{c}}$. If this is the case, then $F^{\et}$ is classically \'etale over $S$.
\end{enumerate}
\end{prop}

As a consequence, to understand objects of $\btlogd$, we first need to understand the extensions of a classical finite \'etale group scheme by a classical finite flat group scheme in the category $\Ab_{\kfl}(S)$ or, equivalently, in $(\fin/S)_{{\mr r}}$ since the latter subcategory is closed by extensions (see \cite[Prop. 2.3]{kat4}).

Let $F',F''\in (\mr{fin}/S)_{\mr{c}}$ and fix a positive integer $n$ that kills both $F'$ and $F''$. We assume $F''$ \'etale and write $F''(1):=F''\otimes_{\Z/n\Z}\Z/n\Z(1)$ where $\Z/n\Z(1)$ denotes the Cartier dual of $\Z/n\Z$.
Let 
\[\EXTc_{S_{\mr{kfl}}}(F'',F')\quad (\text{resp. }\EXTc_{S_{\mr{fl}}}(F'',F'))\]
denote the category of extensions of $F''$ by $F'$ in $(\fs/S)_{\kfl}$ and $(\fs/S)_{\fl}$, respectively. Let
\[\HOMc(F''(1),F')\otimes_{\Z}P^{\mr{gp}}
\]
denote the discrete category associated with the set $\mr{Hom}_S(F''(1),F')\otimes_{\Z}P^{\mr{gp}}$. 
The functor $\Phi_1\colon \EXTc_{S_{\fl}}(F'',F')\to \EXTc_{S_{\kfl}}(F'',F')$, $F^\cl \mapsto F^\cl$, extends to a functor 
\begin{equation}\label{eq1.1}
\Phi=\Phi_\gamma \colon \EXTc_{S_{\fl}}(F'',F')\times \HOMc(F''(1),F')\otimes_{\Z}P^{\mr{gp}}\rightarrow \EXTc_{S_{\kfl}}(F'',F'),
\end{equation}
defined as  \[\Phi(F^{\mr{cl}},\beta):=F^{\mr{cl}}+_{\mr{Baer}}\Phi_2(\beta)\]
where $+_{\mr{Baer}}$ denotes the Baer sum and the functor 
\begin{equation}\label{eq1.2}
\Phi_2= \Phi_{2,\gamma} \colon \HOMc(F''(1),F')\otimes_{\Z}P^{\mr{gp}}\rightarrow \EXTc_{S_{\mr{kfl}}}(F'',F'), \quad \beta\mapsto \Phi_2(\beta),
\end{equation}
is constructed as follows.

For $a\in P^{\mr{gp}}$, let $M_a$ denote the log $1$-motive 
\[ [\Z\stackrel{u_a}{ \longrightarrow }\Gml], \quad u_a(1)=a, \]
where $\Gml$ is Kato's logarithmic multiplicative group on $(\fs/S)_{\kfl}$ (see \cite[Thm. 3.2]{Kato21}).
 %and $u_a(1)=a$. 
 Then $E_{a,n}:=H^{-1}(M_a\otimes^{L}_{\Z}\Z/n\Z)$ fits into a short exact sequence 
\begin{equation}\label{eq.Ea}
0\rightarrow\Z/n\Z(1)\rightarrow E_{a,n}\rightarrow\Z/n\Z\rightarrow0,\end{equation}
which splits Kummer log flat locally. 
Hence tensoring with $F''$ yields another short exact sequence 
\begin{equation}\label{eq.EaF}0\rightarrow F''(1)\rightarrow E_{a,n}\otimes_{\Z/n\Z} F''\rightarrow F''\rightarrow 0.\end{equation}
Now, for any $\Nmono\in\mr{Hom}_S(F''(1),F')$, one defines $\Phi_2(\Nmono\otimes a)$ as the push-out of $E_{a,n}\otimes_{\Z/n\Z} F''$ along $\Nmono$. 
Finally, for any $\beta=\sum_i\Nmono_i\otimes a_i\in \mr{Hom}_S(F''(1),F')\otimes_{\Z}P^{\mr{gp}}$, one defines $\Phi_2(\beta)\in\EXTc_{S_{\mr{kfl}}}(F'',F')$ as the Baer sum of the extensions $\Phi_2(\Nmono_i\otimes a_i)$.

\begin{thm}[Kato]\label{thm1.1}
The functor $\Phi_\gamma $ in \eqref{eq1.1} is an equivalence of categories.
\end{thm}
\begin{proof}
See \cite[Thm. 3.3]{kat4} or \cite[Thm. 3.8]{w-z1}.
\end{proof}

Clearly, the construction of the functor $\Phi_2$ \eqref{eq1.2} involves the chosen chart $\gamma$ and therefore the functor $\Phi$ depends on the chosen chart of $S$. However, once the chart is fixed, $\Phi(F_1^\cl,\beta_1)\simeq \Phi(F_2^\cl, \beta_2)$ if and only if $\beta_1=\beta_2$ and $F_1^\cl\simeq F_2^\cl$. In particular, the following definition makes sense. 

\begin{defn}\label{def.katomonodromyf}
Let $F\in \EXTc_{S_{\kfl}}(F'',F')$. The $\beta \in \HOMc(F''(1),F')\otimes_{\Z}P^{\mr{gp}} $ corresponding to $F$ guaranteed by Theorem \ref{thm1.1} is called the \emph{Kato monodromy} of the extension $F$ of $F''$ by $F'$. 
For $F\in (\fin/S)_{\mr d}$, the \emph{Kato monodromy} of $F$ is defined to be the Kato monodromy of $F$ as an extension of $F^{\et}$ by $F^{\circ}$. If $P^\gp\simeq \Z$, $\beta$ is called the \emph{Kato monodromy map}. 
\end{defn}

We can prove more: once fixed $F$ in $\EXTc_{S_{\kfl}}(F'',F')$, the Kato monodromy of $F$ is essentially independent of the chart chosen on $S$, as explained in the result here below.  
 
\begin{lem}
Let $\gamma'\colon P'_S\to \cM_S$ be another chart. Assume that it induces an isomorphism $\gamma'_{\bar s}\colon P'\xrightarrow{\sim} \cM_{S,\bar s}/\cO_{S,\bar s}^\times$, and set $g=\gamma^{\prime -1}_{\bar s}\circ \gamma_{\bar s} \colon P\xrightarrow{\sim}  P'$. 
Let $\beta'$ denote the Kato monodromy of the extension $F$ in $\EXTc_{S_{\kfl}}(F'',F')$ with respect to $\gamma'$,  and let $g^{\gp}$ be the group envelope of $g$. Then
\[({\rm id}\otimes g^{\gp} ) \colon \HOMc(F''(1),F')\otimes_{\Z} P^{\gp }\rightarrow \HOMc(F''(1),F')\otimes_{\Z} P^{\prime \gp}   \] 
maps the Kato monodromy $\beta$ of $F$ constructed via the chart $\gamma$ to the Kato monodromy $\beta'$ constructed using $\gamma'$.
\end{lem}
\begin{proof}
Note that $\gamma'\circ g_S, \   \gamma\colon P_S\to \cM_S $ are morphisms of sheaves of monoids that induce the same map $P\to \cM_{S, \bar s}  /\cO_{S,\bar s}^\times$. 
Therefore, by Remark \ref{rem.stalks-charts}, $(\gamma'\circ g) - \gamma\colon P\to \cM^\gp_{S,  s}$  factors through $\cO_{S,  s}^\times$ and there exists a $\gamma^\cl\colon P\to \cO_{S,  s}^\times=A ^\times $ such that $\gamma=\gamma'\circ g_S + \gamma^\cl$. 
Let $a\in P^{\gp}$ and $a'=g^\gp (a)$.  Then the short exact sequences $E_{a,n}$  in \eqref{eq.Ea} and the analogous extension $E_{a',n}$  differ by a classical extension over $S$, that is, $E_{a,n}-E_{a',n}\in \EXTc_{S_{\mr{fl}}}(F'',F') $. 
In particular, $\Phi_{2,\gamma}(\Nmono\otimes a)-\Phi_{2,\gamma'}(\Nmono\otimes a')=\Nmono_*(E_{a,n}-E_{a',n})$ is a classical extension for any $\Nmono\in\HOMc(F''(1),F')$. As a consequence, if $\beta=\sum_i \Nmono_i\otimes a_i$ we have 
 \[\Phi_{2,\gamma}\left(\sum_i \Nmono_i\otimes a_i\right)-\Phi_{2,\gamma'}\left(\sum_i \Nmono_i\otimes g^{\gp}(a_i)\right) \in \EXTc_{S_{\mr{fl}}}(F'',F'),\]
 and hence
 \[
F=F^{\cl}+\Phi_{2,\gamma}^n\left(\sum_i \Nmono_i\otimes a_i\right)=F^{\prime \cl}+\Phi_{2,\gamma'}\left(({\rm id}\otimes g^{\gp})\sum_i \Nmono_i\otimes  a_i\right)\]
with $F^{\prime\cl},F^\cl$ suitable extensions in $\EXTc_{S_{\mr{fl}}}(F'',F')$. 
Thus, the Kato monodromy of $F$ with respect to $\gamma'$ is $({\rm id}\otimes g^{\gp})(\sum_i \Nmono_i\otimes  a_i)=({\rm id}\otimes g^{\gp})(\beta)$.
\end{proof}

\begin{cor}\label{cor.katomondromyN}
If $P=P'=\N$, the Kato monodromy map does not depend on the chart.
\end{cor}
\begin{proof}
Clearly, the only possible automorphism of the monoid $\N$ is identity.
\end{proof}

Now, we recall the analogous result for log $p$-divisible groups. 
Let $H'=\varinjlim_{n}H'_n$, $H''=\varinjlim_{n}H''_n$ be two objects in $\btlogc$ (i.e., classical $p$-divisible groups), and assume that $H''$ is \'etale.
Let us denote by                    
\[\EXTc_{S_{\mr{kfl}}}(H'',H')\quad (\text{resp. }\EXTc_{S_{\mr{fl}}}(H'',H'))\] the category of extensions of $H''$ by $H'$ in $\btlogr$ (resp. in $\btlogc$), and by 
\[\HOMc(H''(1),H')\otimes_{\Z}P^{\mr{gp}}
\]
the discrete category associated with the set $\mr{Hom}_S(H''(1),H')\otimes_{\Z}P^{\mr{gp}}$, where $H''(1):=\varinjlim_n H''_n\otimes_{\Z/p^n\Z}\Z/p^n\Z(1)$.

Let 
\[H^{\mr{cl}}= \varinjlim_{n}H_n^{\mr{cl}}\in \EXTc_{S_{\mr{fl}}}(H'',H'),\]
and $\beta\in \HOMc(H''(1),H')\otimes_{\Z}P^{\mr{gp}}$. The element $\beta$ induces a compatible system 
\[\{\beta_n\in \HOMc(H''_n(1),H'_n)\otimes_{\Z}P^{\mr{gp}}\}_{n}.\]
We apply the functor (\ref{eq1.1}) to the pair $(H_n^{\mr{cl}},\beta_n)$ for each $n\geq 1$ and write $\Phi^n_\gamma$ (resp. $\Phi_{2,\gamma}^n$) in place of $\Phi_\gamma$ (resp. $\Phi_{2,\gamma}$) in order to indicate its dependence on $n$. 
Then we get a compatible system  $\{\Phi^n_\gamma(H_n^{\mr{cl}},\beta_n)\}_n$ with 
\[\Phi_\gamma^n(H_n^{\mr{cl}},\beta_n)=H_n^{\mr{cl}}+_{\mr{Baer}}\Phi_{2,\gamma}^n(\beta_n)\in \EXTc_{S_{\mr{kfl}}}(H''_n,H'_n).\]
Note that since $H_n^\cl$ and $E_{a,p^n}\otimes_{\Z/p^n\Z} H_n''$ are both $p^n$-torsion, the same is $\Phi_{2,\gamma}^n(\beta_n)$ and $\Phi_{2,\gamma}(\beta):=\varinjlim_n\Phi_{2,\gamma}^n(\beta_n)$ is an object of $\btlogd$. Therefore 
\[ \varinjlim_n\Phi_\gamma^n(H_n^{\mr{cl}},\beta_n)=\varinjlim_n(H_n^{\mr{cl}}+_{\mr{Baer}}\Phi_{2,\gamma}^n(\beta_n))\]
lies in $\btlogd$. 
We denote $\varinjlim_n\Phi^n_{\gamma}(H_n^{\mr{cl}},\beta_n)$ by $\Phi_{\gamma}(H^\cl,\beta)$. The association of $\Phi_\gamma(H^\cl,\beta)$ to the pair $(H^\cl,\beta)$ gives rise to a functor 
\begin{equation}\label{eq1.3}
\Phi=\Phi_\gamma\colon \EXTc_{S_{\mr{fl}}}(H'',H')\times \HOMc(H''(1),H')\otimes_{\Z}P^{\mr{gp}}\rightarrow \EXTc_{S_{\mr{kfl}}}(H'',H').
\end{equation}

\begin{thm}[Kato]\label{thm1.2}
Let $S$ be as above. 
Assume that there exists a global chart $ \gamma\colon P_S\rightarrow \cM_S$  such that the induced map $P\rightarrow \cM_{S,\bar s}/\mathcal{O}_{S,\bar s}^\times$ is an isomorphism. Let $H',H''\in \btlogc$ with $H''$ \'etale. Then the functor
\[ \Phi_\gamma  \colon \EXTc_{S_{\mr{fl}}}(H'',H')\times \HOMc(H''(1),H')\otimes_{\Z}P^{\mr{gp}}\rightarrow \EXTc_{S_{\mr{kfl}}}(H'',H')\]
in \eqref{eq1.3} is an equivalence of categories.
\end{thm}
\begin{proof}
 This follows from Theorem \ref{thm1.1}.
\end{proof}
As in the finite case, we have a notion of monodromy.

\begin{defn}\label{def.katomonodromy}
Given an object $H$ of $\EXTc_{S_{\kfl}}(H'',H')$, we call the $\beta\in \HOMc(H''(1),H')\otimes_{\Z}P^{\mr{gp}}$ corresponding to $H$ guaranteed by Theorem \ref{thm1.2} the \emph{Kato monodromy} of the extension $H$. If, furthermore, $P^\gp\simeq \Z$ we call it \emph{Kato monodromy map}.
\end{defn}

\subsubsection{Discrete valued base}\label{s.discretevaluedbase} 
%%%%%%%%%%%%%%%%%%%%%%%%%
Now, let $S=\Spec \cO_K$ equipped with the canonical log structure. We fix a uniformizer $\pi$ of $\cO_K$, and thus fix a chart $P:=\N\to \Gamma(S,\mc M_S),1\mapsto\pi$, which satisfies the condition in Theorem \ref{thm1.2}.
For $H\in \btlogd $, let $H^\circ$ (resp.  $H^{\mu}$) be the connected (resp. multiplicative) subgroup of $H$.
As explained in \cite[\S 3.2]{w-z1}, they are classical $p$-divisible groups, and we have a short exact sequence 
\begin{equation}\label{e.connected-etale}
0\rightarrow H^{\circ}\rightarrow H\rightarrow H^{\et}\rightarrow 0 , \end{equation}
with  $H^\et$ classical \'etale.
The \emph{Kato monodromy map} of the log $p$-divisible group $H\in {\btlogd}$ is then defined as the Kato monodromy map
\begin{equation}\label{e.beta}\beta\colon  H^{\et}(1)\rightarrow H^{\circ}
\end{equation}
of $H$ as an extension of $H^{\et}$ by $H^{\circ}$. Since $H^{\et}(1)$ is of multiplicative type, the monodromy $\beta$ of $H$ factors as $H^{\et}(1)\to H^{\mu}\hookrightarrow H^{\circ}$, and, if no confusion arises, we also call the first map Kato monodromy map and denote it by $\beta$ or $\beta^\mu$. Furthermore, by Corollary \ref{cor.katomondromyN}, the Kato monodromy map $\beta$ does not depend on the chart.

\subsection{Kummer log flat cohomology}
%%%%%%%%%%%%%%%%%%%%%%%%%
Let $S$ be an fs log scheme whose underlying scheme is locally noetherian. Let $(\fs/S)_{\fl}$ be the \emph{classical flat site} on $(\fs/S)$, i.e., a covering $\{f_i\colon U_i\to U\}_i$ of an fs log scheme $U$ over $S$ is a set-theoretic covering where the morphisms $f_i$ are strict and their underlying morphisms of schemes are flat and locally of finite presentation \cite[\S4]{Kato21}. We have a forgetful map of sites:
\[\varepsilon \colon (\fs/S)_{\kfl} \to (\fs/S)_{\fl}.\]
In order to understand the cohomology on $(\fs/S)_{\kfl}$, one needs to understand the higher direct images $R^i\varepsilon_*$. 
The following two theorems will be useful for our purpose in this paper. For more results on $R^i\varepsilon_*$ we refer to \cite{zhao-21}.

\begin{thm}\cite[Theorem 4.1]{Kato21}\label{thm-kato}
Let $G$ be a commutative group scheme that is either finite flat or smooth affine over the underlying scheme of $S$.
Then, we have 
\[R^1\varepsilon_{*}G\simeq \varinjlim_n \cHom_S({\Z}/n{\Z}(1), G )\otimes (\Gml/\Gm), \]
where the quotient $\Gml/\Gm$ is taken in $(\fs/S)_{\fl}$.
\end{thm}

\begin{thm}\cite[Theorem 2.3]{zhao-21}\label{thm-zhao}
If $G$ is a torus, then we have 
\begin{enumerate}[(a)]
\item $R^2\varepsilon_{*}G\simeq \varinjlim_n (R^2\epsilon_{*} G)[n]=\bigoplus_{\ell} (R^2\epsilon_{*} G)[\ell^\infty]  $, where $\ell$ varies over all prime numbers; 
\item $(R^2\varepsilon_{*}G)[\ell^r]$ is supported on the locus where the prime $\ell$ is invertible;
\item if $n$ is invertible on $S$, then 
\[ (R^2\varepsilon_{*} G)[n]\simeq G[n](-2)\otimes \wedge^2 ({\mathbb G}_{m,\rm log}/{\mathbb G}_m).\]
\end{enumerate}
\end{thm}

The following example will be used later.

\begin{ex}\label{ex-key} 
Let $R$ be a strictly henselian discrete valuation ring with fraction field $K$ and let $S = \Spec R$ equipped with the canonical log structure. Let $H$ denote both a finite abelian group and the associated constant group scheme over $R$, and let $H^*$ be its Cartier dual. 
For any  resolution  of $H$
\[0\to \Z^r\xrightarrow{\alpha} \Z^r\to H\to 0
\]  
by free abelian groups of finite rank, we get a short exact sequence of group schemes
\[0\to H^*\to \Gm^r\to \Gm^r\to 0.\]
Applying $\varepsilon_*$ to this sequence, we get a long exact sequence
\[\dots \to (R^1\varepsilon_*\Gm)^r\xrightarrow{\beta} (R^1\varepsilon_*\Gm)^r\to R^2\varepsilon_*H^*\to (R^2\varepsilon_*\Gm)^r.\]
By Theorem \ref{thm-kato}, we have $R^1\varepsilon_*\Gm\simeq (\Q/\Z)\otimes_\Z(\Gml/\Gm)$ and the morphism $\beta$ is just $\check{\alpha}\otimes_{\Z}\mr{Id}_{\Gml/\Gm}$, where $\check{\alpha}$ denotes the Pontryagin dual of $\alpha$. Therefore, $\beta$ is surjective. 
Let $(\mr{st}/S)$ be the full subcategory of $(\fs/S)$ consisting of fs log schemes over $S$ whose structure map to $S$ is strict, and denote by $(\mr{st}/S)_{\fl}$ the classical flat site on $(\mr{st}/S)$. 
Then, for any $U\in (\mr{st}/S)$ and any point $u$ of $U$, the stalk of $\Gml/\Gm$ at $\bar{u}$ is either 0 or $\Z$, where $\bar{u}$ is a geometric point above $u$. Thus, the restriction of $\wedge^2 ({\mathbb G}_{m,\rm log}/{\mathbb G}_m)$ to $(\mr{st}/S)_{\fl}$ is zero. By Theorem \ref{thm-zhao}, the restriction of $R^2\varepsilon_*\Gm$ to $(\mr{st}/S)_{\fl}$ is zero.
Therefore, the restriction of $R^2\varepsilon_*H^*$ to $(\mr{st}/S)_{\fl}$ is also zero. 
Then the Leray spectral sequence 
\[E_2^{i,j}=H^i_{\fl}(S,R^j\varepsilon_*H^*)\Rightarrow H^{i+j}_{\kfl}(S,H^*)\]
gives us an exact sequence 
\begin{multline}\label{seven-term exact sequence of Leray spectral sequence}
        0\to  H^1_{\fl}(S,H^*)\to H^1_{\kfl}(S,H^*)\to H^0_{\fl}(S,R^1\varepsilon_*H^*) \\
    \to H^2_{\fl}(S,H^*)\to H^2_{\kfl}(S,H^*)\to H^1_{\fl}(S,R^1\varepsilon_*H^*)    \to  H^3_{\fl}(S,H^*).
\end{multline}
By Theorem \ref{thm-kato}, we have $R^1\varepsilon_*H^*\simeq \check{H}\otimes_{\Z}(\Gml/\Gm)$ with $\check{H}$ the Pontryagin dual of $H$. 
Recall that $R$ is strictly henselian by assumption. 
We have $H^j_{\fl}(S,\Gm)=H^j_{\et}(S,\Gm)=0$ for $j>0$, and thus 
\begin{equation}\label{e.Himult}
H^i_{\fl}(S,H^*)=0 \quad  \text{for} \quad i>1.
\end{equation}
 Then the exact sequence (\ref{seven-term exact sequence of Leray spectral sequence}) gives us an exact sequence
\[0\to H^1_{\fl}(S,H^*)\to H^1_{\kfl}(S,H^*)\to H^0_{\fl}(S,\check{H}\otimes_{\Z}(\Gml/\Gm))\to 0\]
and 
\begin{equation}%\label{e.ex-key}
H^2_{\kfl}(S,H^*)\simeq H^1_{\fl}(S,\check{H}\otimes_{\Z}(\Gml/\Gm))\simeq H^1_{\et}(S,\check{H}\otimes_{\Z}(\Gml/\Gm))=0.
\end{equation}
We are left to compute $H^1_{\kfl}(S,H^*)$. Note that the restriction map 
\[H^1_{\kfl}(S,H^*)\to H^1_\fl(\Spec(K),H^*_K)\] is an isomorphism by \cite[Prop. 3.9]{Gil12}. Over $\Spec(K)$, we have
\[H^i_\fl(\Spec(K),H^*_K)=\begin{cases}0, &\text{ if $i>1$} \\  \check{H}\otimes K^\times, &\text{ if $i=1$} \end{cases}.\]
Since $H^1_{\fl}(S,H^*)\simeq \check{H}\otimes_{\Z}R^\times$, we get an exact sequence 
\begin{equation}\label{e.HHv}
0\to \check{H}\otimes_{\Z}R^\times\to H^1_{\kfl}(S,H^*)\to \check{H}\to 0,
\end{equation}
which can also be induced by the split short exact sequence 
\[0\to R^\times\to K^\times\xrightarrow{\text{valuation}} \Z\to 0.\] 
\end{ex}

\begin{rmk}\label{rmk-key}
In fact, the restriction map $H^1_{\kfl}(S,H^*)\to H^1_\fl(\Spec(K),H^*_K)$ is an isomorphism even when $R$ is merely a discrete valuation ring by \cite[Prop. 3.9]{Gil12}.
\end{rmk}

\subsection{Grothendieck's panach\'ee extensions}\label{s.panach}
%%%%%%%%%%%%%%%%%%%%%%%%%
Let ${\ms C}$ be an abelian category.
We recall here some facts from \cite[Expos\'e IX, Sect. 9.3]{sga72} and \cite[1.5]{Fon87}. 

\begin{defn}\label{def.panachable}
A \emph{panachable sequence} ${\mc S}=(D_1,\cdots, D_5)$ in ${\ms C}$ is a sequence (not exact in general) 
\[D_1\rightarrow D_2\rightarrow D_3\rightarrow D_4\rightarrow D_5\]
such that the induced sequences
\[0\rightarrow D_1\rightarrow D_2\rightarrow D_3\rightarrow 0, \qquad
0\rightarrow D_3\rightarrow D_4\rightarrow D_5\rightarrow 0,
\]
are exact.
\end{defn}

\begin{defn}
    Let $\mc S=(D_1, \cdots, D_5)$ be a panachable sequence in ${\ms C}$. A \emph{panach\'ee extension} of $\mc S$ in ${\ms C}$ is a commutative diagram 
    \[\xymatrix@R=10pt{& & 0\ar[d]& 0\ar[d] & \\ 
0\ar[r] & D_1\ar[r]^j\ar@{=}[d]& D_2\ar[r] \ar[d]& D_3\ar[r]\ar[d]& 0 \\ 0\ar[r]& D_1\ar[r]& D\ar[r]\ar[d]& D_4\ar[r]\ar[d]& 0 \\ 
& & D_5\ar@{=}[r]\ar[d]& D_5\ar[d] \\ & &0&0}\]
in ${\ms C}$ with exact rows and columns.
\end{defn}
The panach\'ee extensions of a panchable sequence $\mc S$ in $\ms C$ form a category (indeed, a groupoid) that we denote as $\EXTPAN_{\ms C}(\mc S)$.  
It follows immediately from the exact sequence
\[
\dots\to \Ext^1_{\ms C}(D_5,D_2)\to \Ext^1_{\ms C}(D_5,D_3)\to \Ext^2_{\ms C}(D_5,D_1)
\]
that $\EXTPAN_{\ms C}(\mc S)$ is not empty if and only if the Yoneda product 
\begin{equation}\label{e.cS}
    c({\mc S}):=D_2\cdot D_4\in \Ext^2_{\ms C}(D_5,D_1)
\end{equation}
is trivial.  Further, the automorphism group of an object $D$ in $\EXTPAN_{\ms C}({\mc S})$ is
\begin{equation}\label{e.aut}
\Hom(D,D)\simeq \Hom_{\ms C}(D_5,D_1) ,
\end{equation} where the left most $\Hom$ means morphisms as panach\'ee extensions.
Finally, by \cite[Expos\'e IX, Prop. 9.3.8]{sga72}, the set $\Extpan_{\ms C}({\mc S})$ of isomorphism classes of objects in $\EXTPAN_{\ms C}({\mc S})$ is a torsor under ${\rm Ext}^1_{\ms C}(D_5,D_1)$ with the action given by 
\begin{align}\label{e.actionpan}
    \omega\colon \Ext^1_{\ms C}(D_5,D_1)\times \Extpan_{\ms C}({\mc S})\to & \ \Extpan_{\ms C}({\mc S}) \\ \nonumber
    ([E],[D])\mapsto &\  j_*[E]+_{\rm Baer}[D]
\end{align}
where the Baer sum is taken as extension classes of $D_5$ by $D_2$.

\section{Obstruction to extending panach\'ee extensions}\label{s.obstruction}
%%%%%%%%%%%%%%%%%%%%%%%%%
The main result of this section is Lemma \ref{lem.pana0} comparing panach\'ee extension in different categories. Let $S:=\Spec \cO_K$ equipped with the canonical log structure. 
Let $(\Sch/K)$ be the category of schemes over $K$, and let $(\Sch/K)_{\fl}$ denote the flat site on $(\Sch/K)$. Let ${\ms C}_{\kfl}$ (resp. ${\ms C}_{\fl}$, resp. ${\ms C}_K$) be the abelian category of sheaves of $\Z/p^n\Z$-modules on the site $S_{\kfl}=(\fs/S)_{\kfl}$ (resp. $S_{\fl}=(\fs/S)_{\fl}$, resp. $(\Sch/K)_{\fl}$). 

Throughout this section $\mc S=(D^\mu,\dots, D^\et)$ denotes a panachable sequence in ${\ms C}_\fl$ (see Def. \ref{def.panachable}) where $D^\mu$
is a finite multiplicative group scheme over $\cO_K$ and $D^\et$ is finite \'etale. 

We have a commutative diagram of continuous functors
 \[\xymatrix{ &S_{\fl}\ar[rd]^{u_3} \\
(\Sch/K)_{\fl}\ar[rr]^{u_1}\ar[ru]^{u_2} &&S_{\kfl}
 }\]
where  $u_3$ is the identity functor and both $u_1
$ and $u_2$ map a $K$-scheme $U$  to $U$ endowed with the trivial log structure. Clearly, $u_1$ and $u_2$ are also cocontinous.
Therefore, the above diagram induces a diagram of topoi 
 \[
 \xymatrix{ &\mr{Sh}(S_{\fl})\ar@<.5ex>[ld]^{f_2} \\
\mr{Sh}((\Sch/K)_{\fl})\ar@<.5ex>[rr]^{g_1}\ar@<.5ex>[ru]^{g_2} &&\mr{Sh}(S_{\kfl})\ar[lu]_{f_3}\ar@<.5ex>[ll]^{f_1}
 },
 \]
 where $f_i=(f_i^{-1}, f_{i*})$ with $f_i^{-1}=(u_i)_s$ exact and $f_{i*}=(u_i)^s $  \cite[tags 00X1, 00XC]{StacksProject} and $g_i=(g_i^{-1}, g_{i*})$ with $g_i^{-1}=(u_i^p)^\sharp$ exact  and $g_{i*}={}_su_{i}$ \cite[tag  00XO]{StacksProject}.  
 By \cite[tag 00XR (1)]{StacksProject}, for any $\mc{F}$ in $\mr{Sh}(S_{\kfl})$ (resp. in $\mr{Sh}(S_{\fl})$) and any $K$-scheme $U$, we have $g_1^{-1}\mc{F}(U)=\mc{F}(U)$ (resp. $g_2^{-1}\mc{F}(U)=\mc{F}(U)$). In particular, if $\mc{F}$ is representable by $Y\in(\fs/S)$, then $g_i^{-1}\mc{F}$ is representable by $Y\times_S\Spec K$ and $f_3^{-1}\mc{F}$ is also representable by $Y$ by \cite[Thm. 3.1]{Kato21}. Furthermore, by \cite[tag 00YV (6)]{StacksProject} $g_1^{-1}$ (resp. $g_2^{-1}$, resp. $f_3^{-1}$) induces an exact functor ${\ms C}_{\kfl}\to {\ms C}_{K}$ (resp. ${\ms C}_{\fl}\to {\ms C}_{K}$, resp. ${\ms C}_{\fl}\to {\ms C}_{\kfl}$).  By construction, we then have a commutative diagram of exact functors 
\begin{equation}\label{restriction.functors}
     \xymatrix{
&{\ms C}_{\fl}\ar[rd]^{f_3^{-1}}\ar[ld]_{g_2^{-1}} \\
{\ms C}_{K} &&{\ms C}_{\kfl}\ar[ll]^{g_1^{-1}}
}.
 \end{equation}
Note that for any two objects $X$ and $Y$ in ${\ms C}_{\kfl}$, if $X$ is \'etale locally represented by a free $\Z/p^n\Z$-module, then we have $\cExt^i_{{\ms C}_{\fl}}(X,Y)=0$ and $\cExt^i_{{\ms C}_{\kfl}}(X,Y)=0$ for any $i>0$. 
By \cite[tag 03FD]{StacksProject} and the local to global Ext spectral sequence in ${\ms C}_{\fl}$ and ${\ms C}_{\kfl}$, we have spectral sequences
\begin{align}\label{local-global-sseqfl}
    H^i_{\fl}(S,\cExt^j_{{\ms C}_{\fl}}(X,Y))&\Rightarrow \Ext_{{\ms C}_{\fl}}^{i+j}(X,Y),\\
 \label{local-global-sseqkfl}
    H^i_{\kfl}(S,\cExt^j_{{\ms C}_{\kfl}}(X,Y))&\Rightarrow \Ext_{{\ms C}_{\kfl}}^{i+j}(X,Y).
\end{align}
Therefore, if $X$ is \'etale locally represented by a free $\Z/p^n\Z$-module, we have
\begin{align}
\label{e.ext-fl}
\Ext_{{\ms C}_{\fl}}^i(X,Y)\simeq &\ H^i_{\fl}(S,\cHom_S(X,Y)), \quad
\\
\label{e.ext-kfl}
\Ext_{{\ms C}_{\kfl}}^i(X,Y)\simeq &\ H^i_{\kfl}(S,\cHom_S(X,Y)), \\
\label{e.ext-kfl2}\Ext^i_{{\ms C}_{K}}(X_K,Y_K)\simeq &\  H^i_{\fl}(\Spec(K) ,\cHom_{S}(X,Y))\end{align}
for any $i\geq0$.

We discuss below results related to the existence of panach\'ee extensions in the above three categories.

\begin{lem}\label{lem-HSLHKL}
Let the notation be as above and set $L=\cHom_S(D^\et,D^\mu)$.
The restriction map $H^1_{\kfl}(S,L)\to H^1_{\fl}(\Spec K,L)$ is an isomorphism.
\end{lem}
\begin{proof}
We prove a bit more, namely that the maps $\theta_2$ and $\theta_5$ in the diagram below are isomorphisms.
 Let $K'$ be a finite unramified Galois field extension of $K$ such that $D^{\et}$ and $(D^{\mu})^*$ become constant over $K'$. Let $\ms{O}_{K'}$ be the ring of integers of $K'$, and endow $S':=\Spec{\ms O}_{K'}$ with the canonical log structure. Then $S'\to S$ is a classical \'etale Galois cover with Galois group $\Gamma:=\mr{Gal}(K'/K)$. 
The \v{C}ech to cohomology spectral sequence \cite[Tag 03OU]{StacksProject} for the cover $S'\to S$ can be expressed via group cohomology as
\[E_2^{i,j}= H^i(\Gamma,H^j_{\kfl}(S',L))\Rightarrow H^{i+j}_{\kfl}(S,L).\]
Similarly, we have a spectral sequence
\[E_2^{i,j}=H^i(\Gamma,H^j_{\fl}(\eta',L))\Rightarrow H^{i+j}_{\fl}(\eta,L),\]
where $\eta':=\Spec K'$ and $\eta:=\Spec K$.
The seven-term exact sequences of the above spectral sequences fit in a diagram that we split into two  diagrams due to lack of space   
\[\xymatrix@R=8pt@C=8pt{
0\ar[r] &H^1(\Gamma,H^0_{\kfl}(S',L))\ar[r]\ar[d]^{\theta_1} &H^1_{\kfl}(S,L)\ar[r]\ar[d]^{\theta_2} &H^0(\Gamma,H^1_{\kfl}(S',L))\ar[r]\ar[d]^{\theta_3} &H^2(\Gamma,H^0_{\kfl}(S',L))\ar[d]^{\theta_4} \\
0\ar[r] &H^1(\Gamma,H^0_{\fl}(\eta',L))\ar[r] &H^1_{\fl}(\eta,L)\ar[r] &H^0(\Gamma,H^1_{\fl}(\eta',L))\ar[r] &H^2(\Gamma,H^0_{\fl}(\eta',L)) \\
\ \\
  &H^2(\Gamma,H^0_{\kfl}(S',L))\ar[d]^{\theta_4}\ar[r] &H^2_{\kfl}(S,L)_{S'}\ar[r]\ar[d]^{\theta_5} &H^1(\Gamma,H^1_{\kfl}(S',L))\ar[r]\ar[d]^{\theta_6} &H^3(\Gamma,H^0_{\kfl}(S',L))\ar[d]^{\theta_7}\\
 &H^2(\Gamma,H^0_{\fl}(\eta',L))\ar[r] &H^2_{\fl}(\eta,L)_{\eta'}\ar[r] &H^1(\Gamma,H^1_{\fl}(\eta',L))\ar[r] &H^3(\Gamma,H^0_{\fl}(\eta',L))\ ; 
}\]
here $H^2_{\kfl}(S,L)_{S'}$ denotes $\ker(H^2_{\kfl}(S,L)\to H^2_{\kfl}(S',L))$ and $H^2_{\fl}(\eta,L)_{\eta'}$ denotes $\ker(H^2_{\fl}(\eta,L)\to H^2_{\fl}(\eta',L))$.  

By \cite[Prop. 3.2.1]{Gil09}, $\theta_3$  and $\theta_6$ are isomorphisms. Since $L$ is finite over $\cO_K$, in particular proper over $\cO_K$, we have $H^0_{\kfl}(S',L)\xrightarrow{\simeq}H^0_{\fl}(\eta',L)$ by the valuative criterion for properness. Therefore $\theta_1$, $\theta_4$  and $\theta_7$ are all isomorphisms. It follows that $\theta_2$ and $\theta_5$ are isomorphisms too, according to the five lemma.  
\end{proof}

The main technical result of this section is the following lemma. It says in particular that   panach\'ee extensions of $\mc S=(D^\mu,\dots, D^\et)$ exist in $\ms C_\kfl$ if and only if they exist for ${\mc S}_K=(D^\mu_K,\dots, D^\et_K)$ in $\ms C_K$. 
 
 \begin{lem}\label{lem.pana0}
 Diagram (\ref{restriction.functors}) induces a commutative diagram of functors
\begin{equation}\label{restriction.functors.panachee}
    \xymatrix{
    &\mr{EXTPAN}_{{\ms C}_{\fl}}({\mc S} )\ar[rd]^{f_3^{-1}}\ar[ld]_{g_2^{-1}} \\
    \mr{EXTPAN}_{{\ms C}_K}({\mc S}_{K })
    &&\ar[ll]^{g_1^{-1}} \mr{EXTPAN}_{{\ms C}_{\kfl}}({\mc S} )
    }
\end{equation} 
Assume that there exists a panach\'ee extension of $\mc S_K$ in ${\ms C}_K$; then the horizontal functor is an equivalence of categories and any panach\'ee extension of ${\mc S}_K$ extends to a unique (up to unique isomorphism) panach\'ee extension of $\mc S$ in ${\ms C}_{\kfl}$.
\end{lem}
\begin{proof}
 The first assertion is immediate. Now assume that there exists a panach\'ee extension $D_K$ of $\mc S_K$ in ${\ms C}_K$. We claim that $\mr{EXTPAN}_{{\ms C}_{\kfl}}({\mc S} )$ is not empty.
 The obstruction to the existence of panach\'ee extensions of ${\mc S}  $ in ${\ms C}_{\kfl} $ is given by the class $c({\mc S} ) \in \Ext_{{\ms C}_{\kfl}}^2(D^\et,D^\mu)\simeq H^2_{\kfl}(S,L)$ with $L=\cHom_S(D^\et,D^\mu)$; see \eqref{e.cS} and \eqref{e.ext-kfl}.

 Claim: There exists a finite Galois unramified field extension $K'$ of $K$ such that $c({\mc S})$ becomes zero in $H^2_{\kfl}(S',L)$, where $S'$ is the spectrum of the ring of integers of $K'$ and we endow it with the canonical log structure. Let $0\to L\to L_1\to L_2\to 0$ be the canonical smooth resolution of $L$, see \cite[Thm. A.5]{Milne1}. Since $H^j_{\fl}(S,L_i)\cong H^j_{\et}(S,L_i)$ for any $j$, we get an exact sequence 
 \[H^1_{\et}(S,L_2)\to H^2_{\fl}(S,L)\to H^2_{\et}(S,L_1).\]
 This exact sequence together with the second part of \eqref{seven-term exact sequence of Leray spectral sequence} and $H^1_{\fl}(S,R^1\varepsilon_*L)\cong H^1_{\et}(S,R^1\varepsilon_*L)$ (see \cite[\href{https://stacks.math.columbia.edu/tag/0DDU}{Tag 0DDU}]{StacksProject}), give us the claim. 

 By the claim, we have $c({\mc S})\in H^2_{\kfl}(S,L)_{S'}$.  If necessary, we enlarge $K'$ so that $D^{\et}$ and $(D^{\mu})^*$ become constant over $K'$ as in the proof of Lemma \ref{lem-HSLHKL}.
Since $c({\mc S}_{K})=0$, by hypothesis, and the map $\theta_5$ in the proof of Lemma \ref{lem-HSLHKL} is an isomorphism, we have $c({\mc S})=0$. We can then fix a panach\'ee extension $D$ of ${\mc S}$ in ${\ms C}_\kfl$.
 
Now, we prove that $g_1^{-1}$ is an equivalence of categories.
The action $\omega$ in \eqref{e.actionpan} gives a bijection 
\[\Ext^1_{{\ms C}_{\kfl}}(D^\et,D^\mu)\isomto \Extpan_{{\ms C}_{\kfl}}({\mc S}),\quad  [E]\mapsto j_*[E]+[D]\]
and similarly, we have a bijection $\Ext^1_{{\ms C}_{K}}(D_K^\et,D_K^\mu)\isomto \Extpan_{{\ms C}_{K}}({\mc S}_{K})$. 
By Remark \ref{rmk-key}, \eqref{e.ext-kfl} and \eqref{e.ext-kfl2} one concludes that  the restriction functor $g_1^{-1}$ induces a bijection 
\begin{equation}\label{e.expanSK} \Extpan_{{\ms C}_{\kfl}}({\mc S})\isomto
\Extpan_{{\ms C}_K}({\mc S}_{K}),\quad j_*[E]+ [D]\mapsto j_*[E_K]+[D_K].
\end{equation}
In particular, the functor $g_1^{-1}$ in \eqref{restriction.functors.panachee} is essentially surjective.

We now prove that $g_1^{-1}$ is fully faithful. First, note that any automorphism of $D_K$ extends uniquely to an automorphism of $D$; in fact, by \eqref{e.aut} and \cite[Chapter III, Lemma 1.1 a)]{Milne1} we have \begin{equation}\label{e.autSK}\Hom(D,D)\simeq \Hom_{S}(D^\et,D^\mu)\simeq \Hom_{\Spec(K)}(D^\et,D^\mu)\simeq \Hom(D_K,D_K).
\end{equation} 
Let $D'$ be a panach\'ee extension of ${\mc S} $ and assume that there exists an isomorphism $ D_K\to D_K'$, i.e., $\Hom(D_K,D_K')$ is not empty. By \eqref{e.expanSK} there exists an isomorphism $\gamma\colon D\to D'$ and hence a commutative diagram 
\[\xymatrix{
\Hom(D,D)\ar[r]^(0.45){\sim}\ar[d]_{\wr}^{\gamma\circ -} &\Hom(D_K,D_K)\ar[d]_{\wr}^{\gamma_K\circ -}  \\
\Hom(D,D')\ar[r] &\Hom(D_K,D'_K)
}
\] where the horizontal arrows are induced by $g_1^{-1}$, the upper one is a bijection by \eqref{e.autSK} and the vertical ones are bijections since $\gamma$ and $\gamma_K$ are isomorphisms. Therefore, the lower one is also bijective and the full faithfulness of $g_1^{-1}$ is clear.
\end{proof}

 We conclude this section with a technical result that will be useful in the study of monodromy pairings.
 
\begin{lem}\label{lem-monodromy-pairing}
Let the notation be as above. Then
   \begin{enumerate}[itemsep=1.5ex,label=(\alph*)]
       \item $\displaystyle{\frac{\Ext^1_{{\ms C}_{\kfl}}(D^{\et},D^{\mu})}{\Ext^1_{{\ms C}_{\fl}}(D^{\et},D^{\mu})}}\simeq \Hom_S(D^\et\otimes_{\Z/p^n}(D^{\mu})^*,\Z/p^n\Z)$.
% \\
% \ 
% \\
% \ 
        \item $\displaystyle{\frac{\Ext^1_{{\ms C}_K}(D_K^{\et},D_K^{\mu})}{\Ext^1_{{\ms C}_{\fl}}(D^{\et},D^{\mu})}}\simeq \Hom_S(D^\et\otimes_{\Z/p^n}(D^{\mu})^*,\Z/p^n\Z)$.
% \\
% \ 
% \\
% \ 
        \item If ${\ms O}_K$ is strictly henselian, then the isomorphism in (b) agrees with that in \cite[\'Exp. IX, Cor. 9.4.4]{sga72}
    \end{enumerate}
\end{lem}
\begin{proof}
Set $L=\cHom_S(D^\et,D^\mu)$.

  (a) We have 
    \begin{eqnarray}\label{e.extfrac}
\frac{\Ext^1_{{\ms C}_{\kfl}}(D^{\et},D^{\mu})}{\Ext^1_{{\ms C}_{\fl}}(D^{\et},D^{\mu})}&\simeq&
\frac{ H^1_{\kfl}(S,L)}{H^1_{\fl}(S,L) } \\
&\simeq &\Hom_S(\mu_{p^n},L)\nonumber \\
&\simeq &\Hom_S(L^*,\Z/p^n\Z)\nonumber \\
&\simeq & \nonumber
 \Hom_S(D^\et\otimes_{\Z/p^n}(D^{\mu})^*,\Z/p^n\Z)
\end{eqnarray}
where the first isomorphism follows by \eqref{e.ext-fl} and \eqref{e.ext-kfl}, the second by \cite[App. D, (D1) and Prop. D.1]{w-z1}.

(b) We have 
    \begin{equation}\label{e.EEHH}
\frac{\Ext^1_{{\ms C}_K}(D_K^{\et},D_K^{\mu})}{\Ext^1_{{\ms C}_{\fl}}(D^{\et},D^{\mu})}  \simeq 
\frac{ H^1_{\fl}(\Spec K,L)}{H^1_{\fl}(S,L) } 
\end{equation}
by \eqref{e.ext-fl} and \eqref{e.ext-kfl2}. Since $H^1_{\kfl}(S,L)\xrightarrow{\sim} H^1_{\fl}( K,L)$ by Lemma \ref{lem-HSLHKL}, the result follows from part (a).    

(c) Since $\ms{O}_K$ is strictly henselian, both $D^{\et}$ and $(D^{\mu})^*$ are constant, finite rank free $\Z/p^n\Z$-modules. Thus, we are reduced to the case that $D^{\et}=(D^{\mu})^*=\Z/p^n\Z$. 
Then we have $L=\mu_{p^n}$. 
By the proof of \cite[Prop. 3.2.1]{Gil09} the isomorphism 
\[H^1_{\kfl}(S,L)/H^1_{\fl}(S,L)\simeq H^1_{\fl}(\Spec K,L)/H^1_{\fl}(S,L)\] induced by  Lemma \ref{lem-HSLHKL} is exactly the first isomorphism in the statement of \cite[\'Exp. IX, Cor. 9.4.4]{sga72}. This finishes the proof.
\end{proof}

%%%%%%%%%%%%%%%%%%%%%%%%%
\section{Log $p$-divisible groups and $p$-divisible groups with sst reduction} \label{s.lpdsst}
%%%%%%%%%%%%%%%%%%%%%%%%%

The main result of this section is Theorem \ref{BZ}, that is part (a) of our Theorem B; we further show that the generic fiber functor respects monodromy (see Theorem \ref{thm.monodromypairings}).   

Let $S:=\Spec \cO_K$ equipped with the canonical log structure. 
Let $H_K$ be a $p$-divisible group over $K$. We say that $H_K$ has \emph{good reduction} if it extends to a $p$-divisible group over ${\ms O}_K$. In this section, we will apply several times Tate's theorem \cite[Thm. 4]{tat1} in case $\mathrm{char}(K)=0$ and de Jong's theorem \cite[Cor. 1.2]{dJong98} in case $\mathrm{char}(K)>0$ stating that the generic fiber functor on $p$-divisible groups is fully faithful.
We have the following definition of semistable reduction of $p$-divisible groups following de Jong \cite[2.2 Definition]{dJong98}.

\begin{defn}\label{def.sstr}
Let $H_K$ be a $p$-divisible group over $K$.  We say that $H_K$ has \emph{semistable reduction} if there exists a filtration $0\subseteq H_K^{\mu}\subseteq H_K^f\subseteq H_K$ such that 
\begin{enumerate}[(a), topsep=0pt]
\item $H_K^f$ (resp. $H_K/H_K^{\mu}$) extends to a $p$-divisible group $\Hone$ (resp. $\Htwo$) over ${\ms O}_K$; 
\item under the condition $(a)$, the morphism $H_K^f\rightarrow H_K\rightarrow H_K/H_K^{\mu}$ extends to a morphism $f\colon \Hone \rightarrow \Htwo$ of $p$-divisible groups over ${\ms O}_K$ with $H^{\mu}:={\rm Ker}(f)$ a multiplicative $p$-divisible group and $H^{\et}:={\rm Coker}(f)$ an \'etale $p$-divisible group. 
\end{enumerate}
\end{defn}
We can depict the data in Definition \ref{def.sstr} as follows
\begin{equation}\label{e.pana}\xymatrix@R=8pt@C=8pt{& & 0\ar[d]& 0\ar[d] & \\ 
0\ar[r] & H_{K}^{\mu}\ar[r]\ar@{=}[d]& H_{K}^{f}\ar[r] \ar[d]& H_{K}^{f}/H_{K}^{\mu}\ar[r]\ar[d]& 0 \\ 0\ar[r]& H_{K}^{\mu}\ar[r]& H_{K}\ar[r]\ar[d]& \Htwo_{K}\ar[r]\ar[d]& 0 \\ & & H_{K}/H_{K}^f\ar@{=}[r]\ar[d]& H_{K}^{\et}\ar[d] \\ & &0&0}\qquad \xymatrix@R=8pt@C=8pt{& &  & 0\ar[d] & \\ 
0\ar[r] & H^{\mu}\ar[r]& \Hone\ar@{.>}[dr]_f\ar[r]  & \Hone/H^{\mu}\ar[r]\ar[d]& 0 \\  &  &  & \Htwo \ar[d]& \\ & &  & H^\et\ar[d] \\ & & &0}
\end{equation} 

The $p$-divisible groups $\Hone, \Htwo$ are denoted as $H_1,H_2$ in \cite{dJong98}.   We have changed the numbering to avoid a conflict of notation with the torsion subgroups of $H$.

\subsection{The canonical filtration}
\label{s.canonicalfiltration}
%%%%%%%%%%%%%%%%%%%%%%%%%

As stated in \cite[2.4 Lemma (i)]{dJong98} any $H_K$ as above admits a \emph{canonical filtration}, i.e., a filtration where $H^\et$ is the \'etale quotient of $\Htwo$ and $H^\mu$ is the multiplicative part of $\Hone$, which is then connected.
Its existence and uniqueness are guaranteed by the lemma below.

\begin{lem}\label{l.filtration}
Let $H_K$ be a $p$-divisible group over $K$ with semistable reduction. Then $H_K$ admits a canonical filtration $0\subseteq H_K^{\mu}\subseteq H_K^f\subseteq H_K$ and for any other filtration $0\subseteq H_K^{\prime \mu}\subseteq H_K^{\prime f}\subseteq H_K$ as in Definition \ref{def.sstr} we have $\Hone$ is isomorphic to the connected component of $H^{\prime \dagger}$. In particular, the canonical filtration is unique.
\end{lem}
\begin{proof}
 We first prove the existence of a canonical filtration by ``extracting'' it from any filtration $0\subseteq H_K^{\mu}\subseteq H_K^f\subseteq H_K$ as in Definition \ref{def.sstr}. 
If $\Hone$ is not connected, we consider the generic fiber of $(\Hone)^{\circ}$ in place of $H_K^f$ and the filtration $0\subseteq H_K^{\mu}\subseteq (\Hone)^{\circ}_K\subseteq H_K$ satisfies the conditions in Definition \ref{def.sstr}.  \begin{equation}\label{e.extraction1}\xymatrix@R=10pt@C=10pt{
&H^\mu \ar[r]&  H^{\dagger, \circ} \ar[dl]\ar@{->>}[r] & (\Hone/H^\mu)^\circ \ar[dl]\ar[d] \\ 
  H^{\mu}\ar[r]\ar@{=}[ru]& \Hone\ \ar@{->>}[r]  &\Hone/H^\mu \ar[d] & \Htwo\ar@{->>}[d]  \\  
 && \Htwo \ar@{=}[ur]\ar@{->>}[d] & H^{\ddagger, \et}\ar@{->>}[dl]\\  
 &&  H^\et   &
 }
 \end{equation}
Note that $(\Hone/H^\mu)^\circ= H^{\dagger, \circ}/H^\mu$. 
If $H^\mu=H^{\dagger,\mu}$, the multiplicative part of $ H^{\dagger, \circ}$, we are done. Otherwise 
let $(\Htwo)'$ be the following push-out of $\Htwo$, \begin{equation}\label{e.extraction2}\xymatrix@R=10pt@C=10pt{&H^\mu \ar[dl ]\ar[r]&  H^{\dagger, \circ} \ar@{=}[dl]\ar@{->>}[r] &  H^{\dagger, \circ}/H^\mu  \ar[dl]\ar[d] \\ 
   H^{\dagger, \mu}\ar[r]&  H^{\dagger, \circ}\ \ar@{->>}[r]  & H^{\dagger, \circ}/ H^{\dagger, \mu} \ar[d] & \Htwo\ar@{->>}[d]\ar[dl]  \\  
 && (\Htwo)'\ar@{->>}[d] &H^{\ddagger, \et}\ar@{=}[dl]  \\
 &&  H^{\ddagger, \et}  &
 }
 \end{equation}
 and note that the filtration $ H^{\dagger, \mu}_K\subseteq  H^{\dagger, \circ}_K\subseteq H_K$ is canonical by construction. The horizontal and vertical sequences in \eqref{e.extraction2} are short exact sequences of $p$-divisible groups, since ${\bf BT}_{{\ms O}_K}$ is closed by extensions \cite[I, (2.4.3)]{Messing}. 

Let now $0\subseteq H_K^{\mu}\subseteq H_K^f\subseteq H_K$ be a canonical filtration and 
$0\subseteq H_K^{\prime \mu}\subseteq H_K^{\prime f}\subseteq H_K$ a filtration as in the statement of the lemma. Then the composition of the inclusion $H_{K}^f\to H_K$ with $H_K\to H^{\prime \et}_K$ is the $0$ map according to Tate's theorem  \cite[Thm. 4]{tat1} in case $\mathrm{char}(K)=0$ and \cite[Cor. 1.2]{dJong98} in case $\mathrm{char}(K)>0$, and hence $H_K^f\subseteq H_K^{\prime f}$.  The inclusion $H_K^f\subseteq H_K^{\prime f}$ corresponds to a map $H^\dagger\to H^{\prime\dagger}$. Since $H^\dagger$ is connected, the map $H^\dagger\to H^{\prime\dagger}$ factors through $(H^{\prime\dagger})^\circ\hookrightarrow H^{\prime\dagger}$. Thus we actually have $H_K^f\subseteq (H^{\prime\dagger})^\circ_K\subseteq H_K^{\prime f}$. Similarly the composition $(H^{\prime\dagger})^\circ_K\to H_K\to H^\et_K$ is the 0 map, and thus $(H^{\prime\dagger})^\circ_K\subseteq H^f_K$. It follows that $(H^{\prime\dagger})^\circ_K= H^f_K$, and thus $H^\dagger\xrightarrow{\simeq}(H^{\prime\dagger})^\circ$ by \cite[Cor. 2]{tat1} and its positive characteristic analogue. 
\end{proof}

%%%%%%%%%%%%%%%%%%%%%%%%%
\subsection{From semistable reduction to log $p$-divisible groups}\label{s.panap}
%%%%%%%%%%%%%%%%%%%%%%%%%

Let $H_K$ be a $p$-divisible group with semistable reduction and fix a filtration $H^\mu_K\subseteq H^f_K\subseteq H_K$ and $p$-divisible groups $\Hone,\Htwo$ as in Definition \ref{def.sstr}.
For any $n$, we have a panachable sequence
\begin{equation}\label{e.Skn}
    {\mc S}_{K,n}=(H_{K,n}^{\mu},\ H_{K,n}^{f},\ H_{K,n}^{f}/H_{K,n}^{\mu},\ H_{K,n}/H_{K,n}^{\mu},\ H_{K,n}/H_{K,n}^{f})
\end{equation} in ${\ms C}_K$ and this extends to a panachable sequence 
\begin{eqnarray}\label{e.Sn}
    {\mc S}^\sst_n :=(H^\mu_n,\ \Hone_n,\ \Hone_n/H^\mu_n,\ \Htwo_n,\ H^\et_n)
\end{eqnarray}
in ${\ms C}_\kfl$ (or in ${\ms C}_{\fl}$). 
Restricting diagrams \eqref{e.pana} to the $p^n$-torsion subgroups, we see that
$H_{K,n}$ has a structure of panach\'ee extension of ${\mc S}_{K,n}$ in ${\ms C}_K$, in particular, the category $\mr{EXTPAN}_{{\ms C}_K}({\mc S}_{K,n})$ is not empty.
Hence, we can apply the results of Section \ref{s.obstruction} and extend $H_K$ to a log $p$-divisible group.

\begin{lem}\label{lem.panaK}
Let $H_K$ be a $p$-divisible group with semistable reduction. The following holds.
\begin{enumerate}[(a)]
\item The  Cartier dual $H_K^*:= \varinjlim_n \cHom (H_{K,n},{\bm \mu}_{p^n})$ of $H_K$ has  semistable reduction.
\item Let ${\mc S}_{K,n}$  and ${\mc S}^\sst_n$ be the panachable sequences attached to a filtration of $H_K$ as in \eqref{e.Skn} and \eqref{e.Sn}.  The panach\'ee extension $H_{K,n}$ of ${\mc S}_{K,n}$ extends to a unique (up to unique isomorphism)  panach\'ee extension $H_n$ of ${\mc S}^\sst_n$ in ${\ms C}_{\kfl}$ and $H_n\in (\mr{fin}/S)_{\mr{d}}$.  
\item For any positive integers $m$ and $n$, the inclusion $H_{K,n}\hookrightarrow H_{K,m+n}$ extends to a unique inclusion $H_n\to H_{m+n}$.
\item $H:=\varinjlim_n H_n $ is an object in $\btlogd$ that extends $H_K$ and is independent of the chosen filtration. 
\item The canonical filtration of $H_K$ is $0\subseteq (H^\mu)_K\subseteq (H^\circ)_K\subseteq H_K$ with $H^\mu$ the multiplicative part of the log $p$-divisible group $H$ in (d) and $H^\circ$ its connected part. 
\end{enumerate}
\end{lem}

\begin{proof}
(a) Applying Cartier duality $(\, )^*$ to both diagrams in \eqref{e.pana}, we see that $H_K^*$ comes equipped with a filtration 
\[0\subseteq H^{*,\mu}_K\subseteq H^{*,f}_K\subseteq H^{*}_K,\] where $H^{*,\mu}_K:=(H^\et)^*_K$ and $H_K^{*,f}:=(\Htwo)^*_K$. Hence, the $p$-divisible group $H^{*}_K$ over $K$ has semistable reduction.

(b) By Lemma \ref{lem.pana0} $H_{K,n}$ extends to a unique (up to unique isomorphism) panach\'ee extension $H_n$ of ${\mc S}^\sst_n$ in ${\ms C}_{\kfl}$.
Note that $H_n$ is an extension of $H_n^\et$ by $\Hone_n$ and hence it lies in $(\mr{fin}/S)_{\mr{r}}$ by \cite[Thm. 9.1]{Kato21}. 
By applying Cartier duality in \eqref{e.pana} and taking the $p^n$-torsion subgroups, we see that the Cartier dual $H_n^*$ of $H_n$ is a panach\'ee extension of  $((H^\et)^*_n, (\Htwo)_n^*,\dots, (\Hone)^*_n , (H^\mu_n)^*)$ and it extends $H_{K,n}^*$ by the uniqueness statement in Lemma \ref{lem.pana0}.
It follows that $H_n\in (\mr{fin}/S)_{\mr{d}}$.

(c) Since $\mr{Ker}(p^n\colon H_{m+n}\to H_{m+n})$ restricts to the panach\'ee extension of $H_{K,n}$, it agrees with $H_n$ by Lemma \ref{lem.pana0}.

(d) It follows from (b) and the proof of (c) that   $H$ satisfies conditions (a) and (c) in Definition \ref{def.logpdivr}. For the surjectivity of the multiplication by $p$, note that $H$ is an extension of the classical $p$-divisible group $H^\et$ over $S$  by $\Hone$.
 Furthermore, $H$ sits in the middle of a diagram \begin{equation*}\xymatrix@R=8pt@C=8pt{  0\ar[r]& H^{\mu}\ar[r] \ar@{=}[d] & \Hone\ \ar[r]\ar[d]  &\Hone/H^\mu\ar[d] \ar[r]&0  \\  
0\ar[r]&H^{\mu}\ar[r]  &H\ar[r]\ar@{->>}[d]& \Htwo \ar@{->>}[d]\ar[r]&0 \\  
 &&H^\et\ar@{=}[r] &  H^\et   &
 }
 \end{equation*} where rows and columns are exact sequences of log $p$-divisible groups; we view it as a compatible system of diagrams on the $p$-power torsion subsheaves in ${\ms C}_{\kfl}$. For the independence statement, it suffices to retrace the proof of Lemma \ref{l.filtration}.  Then diagram \eqref{e.extraction1} and the uniqueness statement in (b) say that $H$ remains the same if we replace the datum of $H^\mu, \Hone,\Htwo, H^\et$ with $H^\mu,  H^{\dagger, \circ},\Htwo, H^{\ddagger, \et}$. On the other hand, diagram \eqref{e.extraction2} and the uniqueness statement in (b) say that $H$ remains the same if we replace the datum of $H^\mu,  H^{\dagger, \circ},\Htwo, H^{\ddagger, \et}$ with $ H^{\dagger,\mu},  H^{\dagger, \circ},(\Htwo)', H^{\ddagger, \et}$. Hence $H$ is isomorphic to the $p$-divisible group in $\btlogd$ constructed using the canonical filtration. 

 The proof of (e) is immediate since \eqref{e.connected-etale} extends the vertical sequence in the middle of the leftmost diagram in \eqref{e.pana} for the canonical filtration.
\end{proof}

\subsection{Monodromy pairings}\label{s.monodromy_pairings}
%%%%%%%%%%%%%%%%%%%%%%%%%

Let $H_K$ be a $p$-divisible group with semistable reduction and $H$ the log $p$-divisible group that extends $H_K$ (see Lemma \ref{lem.panaK}).
Let ${\mc S}_{K,n}=(H^\mu_{K,n}, H^\circ_{K,n}, \dots, H_{K,n}^\et)$ be the panachable sequence \eqref{e.Skn} constructed from the canonical filtration, and let ${\mc S}_n^\sst=(H_n^\mu, H_n^\circ,\dots, H_n^\et) $ be its extension in both ${\ms C}_{\kfl}$ and ${\ms C}_{\fl}$.  
In this section we study three pairings $H_{n}^{\et}\otimes (H_{n}^{\mu})^*\to \Z/p^n\Z$, precisely:
\begin{itemize}
\item the \emph{Grothendieck monodromy pairing}, which measures the obstruction of $H_{K,n}$ lying in the essential image of the functor $g_{2}^{-1}$ in \eqref{restriction.functors.panachee} with $\mc S={\mc S}_n^\sst$;
\item the \emph{logarithmic monodromy pairing}, which measures the obstruction of  $H_n$  lying in the essential image of the functor $f_3^{-1}$ in \eqref{restriction.functors.panachee} with $\mc S={\mc S}_n^\sst$;
\item the \emph{Kato monodromy pairing} induced by the Kato monodromy map in \eqref{e.beta}.
\end{itemize}
We will show that they agree.

The first pairing was constructed in \cite[Expos\'e IX, (9.5.4)]{sga72} via Galois descent from  the strictly henselian case. We present an alternative direct construction below.

Recall that $\pi$ is a fixed uniformizer of ${\ms O}_K$ and the fixed chart on $S$ is $\N_S\to \mc M_S,1\mapsto \pi$. 
 By Theorem \ref{thm1.1},  and Section \ref{s.discretevaluedbase}, the  level $n$-part  of the log $p$-divisible group $H$ is (as element in $\EXT_{S_{\kfl}}(H_n^{\et},H_n^\circ)$)
 \[H_n\simeq \Phi^n(H^\cl_n,\beta_n) :=H^{\cl}_n+_{\mr{Baer}}\Phi_2^n(\beta_n)
 \] for some $H^{\cl}_n\in \EXT_{S_{\fl}}(H_n^{\et},H_n^\circ)$ and some $\beta_n\in\Hom_S(H_n^{\et}(1),H_n^\circ)$, called Kato monodromy map.
 By construction, the extension $\Phi_2^n(\beta_n)$ is the push-out along $\beta_n$ of the sequence \eqref{eq.EaF} with $F''=H_n^\et$.
Since $\beta_n$ factors through a map 
\[\beta_n^\mu\colon \ H_n^\et(1)\to H_n^\mu,\] the extension $\Phi_2^n(\beta_n)$ is also the push-out of an extension 
$H_n^\beta\in \EXT_{S_{\kfl}}
(H_n^\et,H_n^\mu)$.
More precisely, the sheaves  of $\Z/p^n\Z$-modules $H^{\cl}_n$  and $\Phi_2^n(\beta_n)$ fit in the following diagrams
\begin{equation}\label{e.panaclassic}
\xymatrix@R=8pt@C=8pt{0\ar[d] & 0\ar[d] & 0\ar[d]\\
  H_n^\et(1)\ar[d] \ar[r]^{\beta_n^\mu} &H^{\mu}_n\ar[r]^{j_n}\ar[d]& H^{\circ}_n\ar[d] \\ 
    E_{\pi,p^n}\otimes_{\Z/p^n\Z} H^{\et}_n\ar[r]\ar[d] \ar[r] & H^{\beta}_n\ar[r]\ar[d]& \Phi_2^n(\beta_n) \ar[d]  \\ 
  H^\et_n\ar@{=}[r] \ar[d]&H^\et_n\ar@{=}[r] \ar[d]& H^\et_n \ar[d]\\
 0& 0&0}
  \qquad 
  \xymatrix@R=8pt@C=8pt{& & 0\ar[d]& 0\ar[d] & \\ 
0\ar[r] & H^{\mu}_n\ar[r]^{j_n}\ar@{=}[d]& H^{\circ}_n\ar[r]^-{\tau_n} \ar[d]& H^{\circ}_n/H^{\mu}_n\ar[r]\ar[d]& 0 \\ 0\ar[r]& H^{\mu}_n\ar[r]& H_n^\cl \ar[r]\ar[d]& H_{n}^\cl/H^\mu_n \ar[r]\ar[d]& 0 \\ & & H^\et_n \ar@{=}[r]\ar[d]& H^\et_n\ar[d] \\ & &0&0}
\end{equation}
where the left-most vertical sequence is \eqref{eq.EaF} with $F''=H_n^\et$.
As a consequence, the push-out of $\Phi_2^n(\beta_n)$ along $\tau_n$ trivializes, and hence $H_{n}^\cl/H^\mu_n \cong H_{n} /H^\mu_n $ as extensions of $H_n^\et$ by $H^{\circ}_n/H^{\mu}_n$. Therefore,  $H_n^\cl$ is an object of $\EXTPAN_{{\ms C}_{\fl}}({\mc S}_n^\sst)$; in particular, the category $\EXTPAN_{{\ms C}_{\fl}}({\mc S}_n^\sst)$ is not empty.  
    
Let $H_K^{\cl}:=H^{\cl}\times_S\Spec K$. Note that $g_2^{-1}(H^{\cl}_n)$ and $g_1^{-1}(H_n^{\beta})$  are nothing but $H^{\cl}_{K,n}$ and $(H_n^{\beta})_K$ respectively, and we regard them as objects of $\EXTPAN_{{\ms C}_{K}}({\mc S}_{K,n})$ and $\EXT_{{\ms C}_{K}}(H^{\et}_{K,n},H^{\mu}_{K,n})$ respectively. 
By the definition of $\Phi^n$  we have 
    \[ H_{K,n}=\omega_K((H_n^{\beta})_K,H^{\cl}_{K,n}),\]
    where $\omega_K$ is the functor \eqref{e.actionpan} for the category $\EXTPAN_{{\ms C}_{K}}({\mc S}_{K,n})$.

\begin{defn}
    The \emph{Grothendieck monodromy pairing}  of $H_{K,n}$,
    \[c^{\rm Gr}(H_{K,n})\colon \ H_n^\et\otimes_{\Z/p^n}(H_n^{\mu})^*\to \Z/p^n\Z,
    \]
    is the class of $(H_n^{\beta})_K$ in 
    \[\Ext^1_{{\ms C}_K}(H_{K,n}^{\et},H_{K,n}^{\mu})/\Ext^1_{{\ms C}_{\fl}}(H_n^{\et},H_n^{\mu})\simeq \Hom_S(H_n^\et\otimes_{\Z/p^n}(H_n^{\mu})^*,\Z/p^n\Z),\]
    where the isomorphism is given by Lemma \ref{lem-monodromy-pairing}(b). 
\end{defn}
 Note that here we have defined $c^{\rm Gr}(H_{K,n})$ relative to  the panach\'ee extension $H^{\cl}_n\in \EXTPAN_{{\ms C}_{\fl}}({\mc S}^\sst_n)$. 
 As in \cite[\'Exp. IX, p. 108]{sga72} one can define $c^{\rm Gr}(H_{K,n})$ relative to any panach\'ee extension in $\EXTPAN_{{\ms C}_{\fl}}({\mc S}^\sst_n)$; since two different choices differ by an element of $\Ext^1_{{\ms C}_{\fl}}(H_n^{\et},H_n^{\mu})$ under the action $\omega$, the result is independent of the choice. By Lemma \ref{lem-monodromy-pairing}(c), our definition agrees with that of \cite[\'Exp. IX, \S 9.4]{sga72}.
\medskip

 By a similar construction as above, one can define the second pairing.
\begin{defn}
The \emph{logarithmic monodromy pairing} of $H_{n}$ 
\[c^{\log}(H_{n})\colon  H_{n}^{\et}\otimes_{\Z/p^n} (H_{n}^{\mu})^*\to \Z/p^n\Z,\]
is the class of $H_n^{\beta}$ in 
\[\Ext^1_{{\ms C}_{\kfl}}(H_{n}^{\et},H_{n}^{\mu})/\Ext^1_{{\ms C}_{\fl}}(H_n^{\et},H_n^{\mu})\simeq \Hom_S(H_n^\et\otimes_{\Z/p^n}(H_n^{\mu})^*,\Z/p^n\Z).\]
\end{defn}

For constructing the third pairing, consider the map $\beta_n^\mu\colon H^{\et}_n(1)\to H^{\mu}_n$ in \eqref{e.panaclassic} induced by the Kato monodromy map $\beta \colon  H^{\et}(1)\to H^{\circ}$. We have canonical isomorphisms 
\begin{equation}\label{e.dualitiesHom}
    \begin{split}
        \Hom_S(H^{\et}_n(1),H^{\mu}_n)&=\Hom_S(H^{\et}_n,\cHom_S(\mu_{p^n},H^{\mu}_n))\\
&=\Hom_S(H^{\et}_n,\cHom_S( (H^{\mu}_n)^*, \Z/p^n\Z)) \\ &=\Hom_S(H^{\et}_n\otimes_{\Z/p^n}(H^{\mu}_n)^*, \Z/p^n\Z).
    \end{split}
\end{equation}

\begin{defn}
The \emph{Kato monodromy pairing} of $H_{n}$,
\[c(H_{n})\colon  H_{n}^{\et}\otimes_{\Z/p^n} (H_{n}^{\mu})^*\to \Z/p^n\Z,\]
is the pairing associated with $\beta_n^\mu$ via \eqref{e.dualitiesHom}. 
\end{defn}

\begin{lem}
 Let $\delta_{H_{n}^{\beta}}\colon H^{\et}_n\to R^1\varepsilon_*H^{\mu}_n$ be the connecting map for the extension $H^{\beta}_n\in\Ext_{S_{\kfl}}(H^{\et}_n,H^{\mu}_n)$ defined in \eqref{e.panaclassic}. Then we have factorization
    \[\xymatrix{&&\cHom_S(\mu_{p^n},H^{\mu}_n)\ar[d]^\iota \\
    H^{\et}_n\ar[r]_{\delta_{H_{n}^{\beta}}}\ar[rru]^b &R^1\varepsilon_*H^{\mu}_n\ar[r]_-\simeq &\cHom_S(\mu_{p^n},H^{\mu}_n)\otimes_{\Z}\Gmlb
    },\]
    where $b$ is the map corresponding to $\beta_n^\mu$ under the first isomorphism in \eqref{e.dualitiesHom} 
    and $\iota$ is the map $\alpha\mapsto \alpha\otimes[\pi]$. Moreover, $b$ is the only map such that the diagram is commutative.
\end{lem}
\begin{proof}
   See \cite[Lemma 3.3]{w-z1}.
\end{proof}

\begin{thm}\label{thm.monodromypairings}
Let the notation be as above. The Grothendieck monodromy pairing of $H_{K,n}$, the logarithmic monodromy pairing of $H_n$ and the Kato monodromy pairing of $H_n$     agree as pairings $H_n^\et\otimes_{\Z/p^n} (H_n^\mu)^*\to \Z/p^n\Z$.
\end{thm}
\begin{proof}
The first two pairings agree by the very definitions and the commutativity of (\ref{restriction.functors.panachee}).
It remains to show that the logarithmic monodromy pairing and the Kato monodromy pairing of $H_n$  agree. 
As before, we abbreviate $\cHom_S(H_n^\et,H_n^\mu)$ as $L$.   

Let $\varepsilon_{n}\colon \ms{C}_{\kfl}\to \ms{C}_{\fl}$ be the morphism of topoi induced by the forgetful map $\varepsilon\colon (\fs/S)_{\kfl}\to (\fs/S)_{\fl}$ of sites. Let $\mc{F}\in \ms{C}_{\kfl}$. By \cite[tag 03FD, 072W]{StacksProject}, we have 
        \begin{equation}\label{higher-direct-image-is-invariant}
            R^i\varepsilon_* \mc{F}=R^i\varepsilon_{n,*}\mc{F}
        \end{equation}
        as sheaves of abelian groups. We have a spectral sequence 
        \[E_2^{i,j}=\Ext^i_{{\ms C}_{\fl}}(H_{n}^{\et},R^j\varepsilon_{n,*}H_{n}^{\mu})\Rightarrow \Ext^{i+j}_{{\ms C}_{\kfl}}(H_{n}^{\et},H_{n}^{\mu})\]
        by \cite[(2.6)]{zhao-17}, and thus we get another spectral sequence
        \begin{equation}\label{eq.ExtExtss}E_2^{i,j}=\Ext^i_{{\ms C}_{\fl}}(H_{n}^{\et},R^j\varepsilon_{*}H_{n}^{\mu})\Rightarrow \Ext^{i+j}_{{\ms C}_{\kfl}}(H_{n}^{\et},H_{n}^{\mu}).\end{equation}
        Consider the following diagram
        \[\xymatrix@R=12pt@C=8pt{
        &0\ar[d] &0\ar[d] \\
        0\ar[r] &H^1_{\fl}( S,L)\ar[r]^-\simeq\ar[d] &\Ext^1_{{\ms C}_{\fl}}(H_{n}^{\et},H_{n}^{\mu})\ar[r]\ar[d] &H^0( S,\cExt_{{\ms C}_\fl}^1(H_{n}^{\et},H_{n}^{\mu})) \\
        0\ar[r] &H^1_{\kfl}( S,L)\ar[r]^-\simeq_-{\alpha}\ar@{->>}[d]^{\delta} &\Ext^1_{{\ms C}_{\kfl}}(H_{n}^{\et},H_{n}^{\mu})\ar[r]\ar[d]^{\partial} &H^0( S,\cExt_{{\ms C}_{\kfl}}^1(H_{n}^{\et},H_{n}^{\mu})) \\
        &H^0(S,R^1\varepsilon_*L) &\Hom_S(H_{n}^{\et},R^1\varepsilon_{*}H_{n}^{\mu})
        }\]
        with exact rows and columns, where the two rows are the three-term exact sequence of the spectral sequences \eqref{local-global-sseqfl} and \eqref{local-global-sseqkfl}, the first column is the three-term exact sequence of the Leray spectral sequence, and the second column is the three-term exact sequence of \eqref{eq.ExtExtss}.
        One can check that the left upper square is commutative. The map $\delta$ is surjective by \cite[App. D, Prop. D.1 (2)]{w-z1}. We identify $H^0(S,R^1\varepsilon_*L)$ with 
        \[\Hom_S(\mu_{p^n},L)=\Hom_S(H_n^\et\otimes_{\Z/p^n}(H_n^{\mu})^*,\Z/p^n\Z)\]
        as in \cite[App. D, Prop. D.1 (1)]{w-z1}, and the latter is clearly finite. Then the element $\delta(\alpha^{-1}(H^{\beta}_n))$ is exactly the logarithmic monodromy pairing by definition. Since 
        \begin{align*}
            \Hom_S(H_{n}^{\et},R^1\varepsilon_{*}H_{n}^{\mu})=&\Hom_S(H_{n}^{\et},\Hom_S(\mu_{p^n},H_{n}^{\mu})\otimes\Gmlb) \\
            \simeq&\Hom_S(H_{n}^{\et},\Hom_S(\mu_{p^n},H_{n}^{\mu}))\\
            \simeq&\Hom_S(H_{n}^{\et}(1),H_{n}^{\mu}) \\
            \simeq&\Hom_S(H_n^\et\otimes_{\Z/p^n}(H_n^{\mu})^*,\Z/p^n\Z),
        \end{align*}
        where the second isomorphism follows from \cite[Lem. 3.7]{w-z1}, the groups $H^0(S,R^1\varepsilon_*L)$ and $\Hom_S(H_{n}^{\et},R^1\varepsilon_{*}H_{n}^{\mu})$ are finite of the same order.
         Then the above commutative diagram implies that $\partial$ is also surjective.
        We identify $\Hom_S(H_{n}^{\et},R^1\varepsilon_{*}H_{n}^{\mu})$ with $\Hom_S(H_{n}^{\et}(1),H_{n}^{\mu})$, then $\partial(H^{\beta}_n)$ is just the Kato monodromy map by definition. It follows that under the isomorphism 
        \[\Hom_S(H_n^\et\otimes_{\Z/p^n}(H_n^{\mu})^*,\Z/p^n\Z)\xrightarrow{\simeq}\Hom_S(H_{n}^{\et}(1),H_{n}^{\mu})\]
        induced by the isomorphism $\alpha$, the logarithmic monodromy pairing is mapped to the Kato monodromy map. Therefore, the logarithmic monodromy pairing agrees with the Kato monodromy pairing. 
\end{proof} 

\begin{defn}\label{p-adic Kato pairing}
    The \emph{Kato monodromy pairing} of $H$
 \[c(H)\colon T_p(H^{\et})\otimes_{\Z_p} T_p((H^{\mu})^*) \to \Z_p\]
 is now defined by passing to the inverse limit on $c(H_n)$  or, equivalently, on $c^{\log}(H_n)$. Similarly one defines the Grothendieck monodromy pairing $c^{\rm Gr}(H_K)$ by passing to the inverse limit on $c^{\rm Gr}(H_{K,n})$.   
\end{defn}

Note that by Theorem \ref{thm.monodromypairings}  $c(H)$ and $c^{\rm Gr}(H_K)$ are the same pairing.

\subsection{Criterion for semistable reduction of $p$-divisible groups} 
%%%%%%%%%%%%%%%%%%%%%%%%%
\begin{thm}\label{BZ}
There is an equivalence of categories 
 \[\btlogd\to \bt_{K}^\sst,  \quad H\mapsto H\times_S\Spec K.\] 
\end{thm}
\begin{proof} 
 Let $H\in \btlogd$. Then $H_n\to S$ is Kummer. Since the log structure of $S$ is supported on the closed point, the map $H_n\times_S\Spec K\to \Spec K$ has to be strict. It follows that $H\times_S\Spec K$ is a classical $p$-divisible group over $K$ and hence, we have a functor 
 \[(\ )_K\colon \btlogd\to \bt_K,\quad H\mapsto H\times_S\Spec K .\]
 It remains to prove that this functor is fully faithful with essential image  $\bt_{K}^\sst$.

  We determine the image of the functor. In the connected-\'etale decomposition \eqref{e.connected-etale} of $H$ both $H^\circ$ and $H^{\et}$ lie in $\btlogc$, and $H^\circ$ (resp. $H^{\et}$) is connected (resp. \'etale) (see \cite[Prop. 2.7 (3)]{kat4} or \cite[Prop. 3.9]{w-z1}). Let $H^{\mu}\subset H^\circ$ be the maximal multiplicative $p$-divisible subgroup of $H^\circ$. 
Now the connected log $p$-divisible group $(H/H^\mu)^\circ$ over $S$ has trivial multiplicative subgroup; hence, the monodromy morphism associated with $H/H^\mu$ in Kato's classification of log $p$-divisible groups must be trivial (see \cite[Cor. 3.14]{w-z1}). In particular, $H/H^\mu$ is classical.
Then the $p$-divisible group $H_K:=H\times_S\Spec K$ has a filtration 
$0\subseteq H_K^{\mu}\subseteq H_K^f\subseteq H_K$
which verifies the condition of semistable reduction with $\Hone=H^\circ$ and $\Htwo=H/H^\mu$. 
Therefore the image of the functor $(\ )_K$ is contained in $\mathbf{BT}_{K}^{\mathrm{st}}$, and this is the essential image by Lemma \ref{lem.panaK}(d).

  We show that the functor is faithful. Let $f\colon G\to H$ be a morphism in $\btlogd$ such that $f_K:=f\times_S\Spec K=0$. It suffices to show that $f=0$. Let $0\to G^\circ\to G\to G^\et\to0$ be the connected-\'etale decomposition of $G$. Then $f$ induces morphisms $f^\circ\colon G^\circ\to H^\circ$ and $f^\et\colon G^\et\to H^\et$. Since $f^\circ\times_S\Spec K=0$ and $f^\et\times_S\Spec K=0$, we must have $f^\circ=0$ and $f^\et=0$ by \cite[Thm. 4]{tat1} and \cite[Cor. 1.2]{dJong98}. 
  It follows that $f$ factors as $G\to G^\et\xrightarrow{\bar{f}}H^\circ\to H$. The vanishing of $f_K$ implies that $\bar{f}_K=0$. Applying \cite[Thm. 4]{tat1} and \cite[Cor. 1.2]{dJong98} again, we get $\bar{f}=0$. Therefore $f=0$.

  At last, we show that the functor is full. Let $G$ and $H$ be in $\btlogd$, and let $g_K\colon G_K\to H_K$ be a morphism of $\bt_K$. It is enough to extend $g_K$ into a morphism of $\btlogd$. The composition $\gamma_K\colon G^\circ_K\to G_K\xrightarrow{g_K} H_K\to H^\et_K$ extends to a morphism $\gamma\colon G^\circ\to H^\et$ by \cite[Thm. 4]{tat1} and \cite[Cor. 1.2]{dJong98}. Since $G^\circ$ is connected and $H^\et$ is \'etale, $\gamma$ has to be trivial and thus $\gamma_K=0$. Then one can see that $g_K$ induces $G^\circ_K\to H^\circ_K$ and $G^\et_K\to H^\et_K$ which extend to morphisms $g^\circ\colon G^\circ\to H^\circ$ and $g^\et\colon G^\et\to H^\et$ respectively. In order to extend $g_K$, it suffices to identify the push-forward of $G$ as an extension of $G^\et$ by $G^\circ$ along $g^\circ$ and the pull-back of $H$ as an extension of $H^\et$ by $H^\circ$ along $g^\et$, as depicted below
\[\xymatrix{
0\ar[r] &G^\circ\ar[r]\ar[d]_{g^\circ} &G\ar[r]\ar[d] &G^\et\ar[r]\ar@{=}[d] &0\\
0\ar[r] &H^\circ\ar[r]\ar@{=}[d] &(g^\circ)_*G\ar[r]\ar@{..>}[d]^{?}_{\cong} &G^\et\ar[r]\ar@{=}[d] &0\\
0\ar[r] &H^\circ\ar[r]\ar@{=}[d] &(g^\et)^*H\ar[r]\ar[d] &G^\et\ar[r]\ar[d]^{g^\et} &0\\
0\ar[r] &H^\circ\ar[r] &H\ar[r] &H^\et\ar[r] &0\\
}.
\]
Note that the dotted arrow exists over $K$ and is an isomorphism. So we are reduced to the case that $G\in\btlogd$ such that $G^\circ=H^\circ$, $G^\et=H^\et$, $g^\circ=1_{H^\circ}$, $g^\et=1_{H^\et}$, and $g_K$ is an isomorphism.
Let $H^\mu$ be the multiplicative part of $H^\circ$. The similar argument in the beginning of this part shows that $g_K$ induces $1_{H^\mu}$ and an isomorphism $g_{K,\mu}\colon G_K/H^\mu_K\to H_K/H^\mu_K$, and the latter extends to a unique isomorphism $g_{\mu}\colon G/H^\mu\to H/H^\mu$ by \cite[Thm. 4]{tat1} and \cite[Cor. 1.2]{dJong98}.
We identify $G/H^\mu$ with $H/H^\mu$ through $g_\mu$, and denote it by $H^\ddagger$. For each positive integer $n$, let ${\mc S}:= (H^\mu_{n},H^\circ_{n},H^\circ_{n}/H^\mu_{n},H_{n}^\ddagger,H^\et_n)$. Then $G_n$ and $H_n$ are objects of $\mr{EXTPAN}_{{\ms C}_{\kfl}}({\mc S} )$ that restrict to isomorphic objects in $\mr{EXTPAN}_{{\ms C}_K}({\mc S}_K )$. Therefore, by Lemma \ref{lem.pana0} there is a unique isomorphism $g_n\colon G_n\simeq H_n$ that extends $g_{K,n}$, and thus $g_K$ extends to an isomorphism $g\colon G\to H$.
\end{proof}

%%%%%%%%%%%%%%%%%%%%%%%%% 
\section{Fontaine's conjecture for log $p$-divisible groups}
%%%%%%%%%%%%%%%%%%%%%%%%%

Let $\cO_K$ satisfy the stronger assumption $(\ast)$ from the Introduction and let $S$ be $\Spec \cO_K$ equipped with the canonical log structure. 
 
 In this section, we will prove the second part of Theorem \ref{main}, which is the logarithmic analogue of Fontaine's conjecture on Galois representations associated with $p$-divisible groups.

\subsection{From log $p$-divisible groups to Galois representations}
%%%%%%%%%%%%%%%%%%%%%%%%%

Now for a log $p$-divisible group $H\in \btlogd$, we denote by $V_{p}(H)\in {\bf Rep}_{\Q_p}(\Galk) $ the $p$-adic Galois representation attached to the generic fiber $H_K$ of $H$.

Let $0\to H^\circ\to H\to H^{\et}\to 0$ be the connected-\'etale decomposition of $H\in \btlogd$. By Theorem \ref{thm1.2}, there exists a classical $p$-divisible group $H^{\cl}$ over $S$ that is an extension of $H^{\et}$ by $H^\circ$, and a homomorphism $\beta\colon H^{\et}(1)\to H^\circ$, such that $H=\Phi(H^{\cl},\beta)$. By the construction of $\Phi(H^{\cl},\beta)$, the $p$-adic Galois representation $V_{p}(H)$ can be constructed from the pair $(V_{p}(H^\cl),\beta)$ as follows. 

Let 
\begin{equation}\label{the rep. associated to the log 1-mot M_pi}
    V_\pi:=V_{p}(M_{\pi}[p^\infty])\in {\bf Rep}_{\Q_p}(\Galk)
\end{equation}
be the $p$-adic Galois representation associated to the log $p$-divisible group $M_{\pi}[p^\infty]$ of the log $1$-motive $M_\pi=[\Z\stackrel{1\mapsto \pi}{\longrightarrow} \Gml]$. Note that $V_{\pi}$ is also the $p$-adic Galois representation associated to the Tate curve with $q$-invariant $\pi$. 
Thus $V_{\pi}$ is a semistable representation by \cite[IV.5.4]{berg1}  and we have a short exact sequence of $p$-adic Galois representations
\begin{equation}\label{eq1}
0\to\Q_p(1)\to V_\pi\to \Q_p\to 0.
\end{equation}
Tensoring with the unramified representation $V_{p}(H^\et)$ of $\Galk$, we get a short exact sequence
\begin{equation}\label{eq2}
0\to V_{p}(H^\et)(1)\to V_{\pi}\otimes_{\Q_p}V_{p}(H^\et)\to V_{p}(H^\et)\to 0.
\end{equation}

The $p$-adic Galois representation $V_{p}(\Phi_2(\beta))$ associated to the log $p$-divisible group $\Phi_2(\beta)$ is given by the push-out
\begin{equation}\label{eq3}
\xymatrix{
0\ar[r] &V_{p}(H^\et)(1)\ar[r]\ar[d]^{V_{p}(\beta)} &V_{\pi}\otimes_{\Q_p}V_{p}(H^\et)\ar[r]\ar[d] &V_{p}(H^\et)\ar[r]\ar@{=}[d] &0  \\
0\ar[r] &V_{p}(H^\circ)\ar[r] &V_{p}(\Phi_2(\beta))\ar[r] &V_{p}(H^\et)\ar[r] &0 
}.
\end{equation}
Then $V_{p}(H)$ is given by the Baer sum of $V_{p}(\Phi_2(\beta))$ and $V_{p}(H^{\cl})$ as extensions of $V_{p}(H^{\et})$ by $V_{p}(H^\circ)$.

Aiming to prove that $V_{p}(H)$ is semistable, we first check that usual operations on extensions of $p$-adic Galois representations respect semistability.

\begin{lem}\label{strep1}
Let $0\to V_1\to V_2\to V_3\to 0$ be an extension of $p$-adic Galois representations with $V_1$ and $V_3$ semistable. 
\begin{enumerate}[(a)]
\item $V_2$ is semistable if and only if the sequence 
\[0\to D_{\rm st}(V_1)\to D_{\rm st}(V_2)\to D_{\rm st}(V_3)\to 0\]
is exact, where $D_{\sst}(V):=(V\otimes_{\Q_p} B_{\sst})^{\Galk}$ is the functor from \cite[\S4.1]{CF00}, see also \eqref{functor.Dst}.
\item Let $f\colon V_1\to V_1'$ be a homomorphism of semistable $p$-adic Galois representations. Assume that $V_2$ is semistable. Then the push-out $f_*V_2$ of $V_2$ along $f$ is also semistable.
\item Let $g\colon V_3'\to V_3$ be a homomorphism of semistable $p$-adic Galois representations. Assume that $V_2$ is semistable. Then the push-back $g^*V_2$ of $V_2$ along $g$ is also semistable.
\item Let $0\to V_1\to V'_2\to V_3\to 0$ be another extension of $p$-adic Galois representations. Assume that both $V_2$ and $V_2'$ are semistable. Then the Baer sum of $V_2$ and $V_2'$ as extensions is also semistable.
\end{enumerate}
\end{lem}
\begin{proof}
(a) The sequence is obviously left exact. It is exact if and only if 
\[\mr{dim}_{K_0}(D_{\sst}(V_2))=\mr{dim}_{K_0}(D_{\sst}(V_1))+\mr{dim}_{K_0}(D_{\sst}(V_3)).\]
Since $V_1$ and $V_3$ are semistable, we get 
\[\mr{dim}_{K_0}(D_{\sst}(V_1))+\mr{dim}_{K_0}(D_{\sst}(V_3))=\mr{dim}_{\Q_p}(V_1)+\mr{dim}_{\Q_p}(V_3)=\mr{dim}_{\Q_p}(V_2).\]
Since $V_2$ is semistable if and only if $\mr{dim}_{K_0}(D_{\sst}(V_2))=\mr{dim}_{\Q_p}(V_2)$, the result follows.

(b) Since $f_*V_2$ is a push-out, we have the following commutative diagram
\[\xymatrix{
0\ar[r] &D_{\sst}(V_1)\ar[r]\ar[d]^{D_{\sst}(f)} &D_{\sst}(V_2)\ar[r]\ar[d] &D_{\sst}(V_3)\ar@{=}[d]  \\
0\ar[r] &D_{\sst}(V_1')\ar[r] &D_{\sst}(f_*V_2)\ar[r] &D_{\sst}(V_3) 
}\]
with exact rows. Since $V_2$ is semistable, the last map of the upper row of the diagram is actually surjective by (a). And thus so is the last map of the lower row. Again by (a), $f_*V_2$ is semistable.

(c) Since $g^*V_2$ is a pull-back, we have the following commutative diagram
\[\xymatrix{
0\ar[r] &D_{\sst}(V_1)\ar[r]\ar@{=}[d] &D_{\sst}(g^*V_2)\ar[r]\ar[d] &D_{\sst}(V_3')\ar[d]^{D_{\sst}(g)}\ar[r]^-{\delta'} &H^1(\Galk,V_1\otimes_{\Q_p}B_{\sst})\ar@{=}[d]\\
0\ar[r] &D_{\sst}(V_1)\ar[r] &D_{\sst}(V_2)\ar[r] &D_{\sst}(V_3)\ar[r]^-{\delta} &H^1(\Galk,V_1\otimes_{\Q_p}B_{\sst})
}\]
with exact rows. Since $V_2$ is semistable, we get $\delta=0$ by (a). Therefore $\delta'=0$. Again by (a), we have that $g^*V_2$ is semistable.

(d) Since the Baer sum is constructed out of product, push-out and pull-back, the result follows from (a), (b) and (c).
\end{proof}

\begin{lem}\label{strep2}
The Galois representation $V_{\pi}\otimes_{\Q_p}V_{p}(H^\et)$ is semistable.
\end{lem}
\begin{proof}
Since $H^\et$ is a classical $p$-divisible group over $S$, the representation $V_{p}(H^\et)$ is crystalline, in particular semistable. By \cite[Prop. 4.2]{CF00}, the tensor product $V_{\pi}\otimes_{\Q_p}V_{p}(H^\et)$ of two semistable representations is also semistable.
\end{proof}

\begin{prop}\label{the rep. associated to a log BT is sst with HT weight 0 and 1}
Let $H\in \btlogd$. Then $V_{p}(H)\in {\bf Rep}_{\Q_p}(\Galk)$ is semistable with Hodge-Tate weights in $\{0,1\}$.
\end{prop}
\begin{proof}
Since $H^{\cl}$ is a classical $p$-divisible group, $V_{p}(H^{\cl})$ is semistable  (see Theorem \ref{classicalA}). As $V_{p}(H)$ is given by the Baer sum of $V_{p}(\Phi_2(\beta))$ and $V_{p}(H^{\cl})$ as extensions of $V_{p}(H^{\et})$ by $V_{p}(H^\circ)$, we are reduced to showing that $V_{p}(\Phi_2(\beta))$ is semistable by Lemma \ref{strep1} (d). But the semi-stability of $V_{p}(\Phi_2(\beta))$ follows from Lemma \ref{strep2} and Lemma \ref{strep1} (b).

 Since $V_{p}(H^\et)$ and $V_{p}(H^\circ)$ are $p$-adic Galois representations associated to classical $p$-divisible groups, they are crystalline representations with Hodge-Tate weights in $\{0,1\}$ by Theorem \ref{classicalA}. The representations $V_{p}(H^\circ)$, $V_{p}(H^\et)$ and $V_{p}(H)$ are all Hodge-Tate. By \cite[Prop. 1.6 (iii)]{Fon82}, the functor $D_{\mathrm{HT}}(-):=(-\otimes_{\Q_p}B_{\mathrm{HT}})^{\Galk}$ is an exact functor on the category of Hodge-Tate representations, therefore we have a short exact sequence 
 \[0\to D_{\mathrm{HT}}(V_{p}(H^\circ))\to D_{\mathrm{HT}}(V_{p}(H))\to D_{\mathrm{HT}}(V_{p}(H^\et))\to 0.\]
 It follows that the  Hodge-Tate weights of $V_{p}(H)$ are in $\{0,1\}$.
\end{proof}

\subsection{From Galois representations to log $p$-divisible groups}
%%%%%%%%%%%%%%%%%%%%%%%%%

In this subsection we  associate to any $\rho\in \Repz^{\sst, \{0,1\}}(\Galk)$ a logarithmic $p$-divisible group in $\btlogd$.  
The key ingredient is Fargues' theory of $p$-divisible rigid analytic groups in \cite{Far19, Far22}. In this subsection, analytic space means paracompact strictly $K$-analytic space in the sense of Berkovich, or equivalently, quasi-separated rigid $K$-analytic space that has an admissible affinoid covering of finite type. The equivalence is locally described by associating Berkovich spectrum $\mathcal M(A)$ with the maximal spectrum ${\rm Sp}(A)$ when $A$ is a strictly $K$-affinoid algebra.

Recall that a $p$-divisible rigid analytic $K$-group is a commutative rigid analytic $K$-group such that the $p$-multiplication is topologically nilpotent, finite, and surjective \cite[Def. 1.1]{Far22}. 
One should not confuse $p$-divisible rigid analytic groups with $p$-divisible groups; indeed any object in $\bt_K$ produces a $p$-divisible rigid analytic $K$-group, but the converse does not hold: $\Ga^\rig$ is a counterexample.

Let $\bt_K^\rig$ be the category of $p$-divisible rigid analytic $K$-groups. This notation is taken from \cite[\S1]{Far22}, and the corresponding notation in \cite[\S2.1, Def. 3]{Far19} is $\mathcal{R}_K$. 
Let $C$ be the completion of a fixed algebraic closure of $K$. By \cite[\S1]{Far22} or more precisely \cite[Cor. 17]{Far19}, there is an equivalence of categories between $\bt_K^\rig$ and the category of triples $(\Lambda,W,\alpha)$, where $\Lambda$ is a (continuous) representation of $\Galk$ on a finite rank free $\Z_p$-module, $W$ is a finite dimensional $K$-vector space, and 
\[\alpha\colon W_C(1):=W\otimes_KC(1)\to\Lambda\otimes_{\Z_p}C=:\Lambda_{C}\]
is a $C$-linear map which is compatible with the Galois actions. Given a triple $(\Lambda,W,\alpha)$, let $G^\rig$ be the corresponding $p$-divisible rigid analytic $K$-group. Then as stated in \cite[\S1]{Far22} $\Lambda$ can be recovered from $G^\rig$ as 
\begin{equation}\label{recover representation from associated p-div rig ana gp}    \Lambda=T_p(G^\rig[p^\infty]),
\end{equation}
where $G^{\rig}[p^n]$ is the $p^n$-torsion subgroup of $G^{\rig}$ and $G^{\rig}[p^\infty]:=\varinjlim_nG^{\rig}[p^n]$. Note that $G^\rig[p^\infty]$ can be viewed as an object in $\bt_K$ by \cite[\S 1.2, Cor. 2]{Far19}.

Let $\rho\colon \Galk\to\mathrm{GL}(\Lambda)$ be in $\Repz^{\sst,\{0,1\}}(\Galk)$, and consider the triple 
\[(\Lambda,(\Lambda_C(-1))^{\Galk},\alpha)\]
with $\alpha$ the canonical inclusion.
Let $G^{\rig}(\rho)$ be the $p$-divisible rigid analytic $K$-group corresponding to $(\Lambda,(\Lambda_C(-1))^{\Galk},\alpha)$, see the beginning of \cite[\S 2]{Far22}; 
here we add a superscript ${}^\rig$ to stress that it is a rigid analytic object.
We regard $G^{\rig}(\rho)[p^\infty]$ as an object in $\bt_K$, and $\rho$ is just $T_p(G^{\rig}(\rho)[p^\infty])$ by \eqref{recover representation from associated p-div rig ana gp}.
We are going to show that $G^{\rig}(\rho)[p^\infty]$ has semistable reduction and thus extends to a unique $G\in \btlogd$ by Theorem \ref{BZ}.

\begin{prop}\label{the $p$-div gp associated to a sst rep. with HT weights 0 and 1 is sst}
 We have $G^{\rig}(\rho)[p^\infty]\in {\bf BT}_{K}^{\mr{st}}$.
\end{prop}
\begin{proof}
By \cite[\S 1.4, Prop. 8]{Far19}, $G^{\rig}(\rho)$ fits into a short exact sequence
\[0\to G^{\rig}(\rho)^\circ\to G^{\rig}(\rho)\to \underline{\pi}_0(G^{\rig}(\rho))\to 0\]
of sheaves of abelian groups on the big \'etale site of $\mathcal M(K)$, where $G^{\rig}(\rho)^\circ$ denotes the identity component of $G^{\rig}(\rho)$ and $\underline{\pi}_0(G^{\rig}(\rho))$ is the \'etale analytic group of connected components. By \cite[\S 2.8, Prop. 21]{Far19}, the above short exact sequence gives rise to a short exact sequence
\[0\to G^{\rig}(\rho)^\circ[p^\infty]\to G^{\rig}(\rho)[p^\infty]\to \underline{\pi}_0(G^{\rig}(\rho))\to 0\]
of $p$-divisible rigid analytic groups, which can also be regarded as a short exact sequence of $p$-divisible groups over $K$ by \cite[\S 2.10, Cor. 13]{Far19}.

By \cite[Prop. 2.1]{Far22},
since $\rho$ is semistable with Hodge-Tate weights in $\{0,1\}$, the action of $\Galk$  on $\underline{\pi}_0(G^{\rig}(\rho))(\overline{K})$ is unramified and $G^{\rig}(\rho)^\circ$ is isomorphic to the open unit ball $\mathring{\mathbf{B}}^d_K$  of dimension $d$. Therefore, the $p$-divisible group  $\underline{\pi}_0(G^{\rig}(\rho))$ in $\bt_K$ extends to an \'etale $p$-divisible group $G^\et$ over $\cO_K$. 
By \cite[\S 6, Thm. 6.1]{Far19}, there exists a $p$-divisible formal group $F$ over $\cO_K$ whose associated rigid analytic group is $G^{\rig}(\rho)^\circ$. 
Let $G^\circ:=F[p^\infty]$ be the (formal) $p$-divisible group over $\cO_K$ associated with $F$ (see \cite[Introduction]{Far19}), and let $G^\mu$ be the multiplicative part of $G^\circ$. Then we have $G^\circ_K=G^{\rig}(\rho)^\circ[p^\infty]$, $G^\et_K=\underline{\pi}_0(G^{\rig}(\rho))$, and a commutative diagram
\begin{equation}\label{pan-ext diagram associated a sst represenation of HT weights 0 and 1}
\xymatrix@R=10pt{& & 0\ar[d]& 0\ar[d] & \\ 
0\ar[r] & G^\mu_K\ar[r]\ar@{=}[d]& G^\circ_K\ar[r] \ar[d] \ar[rd] & G^\circ_K/G^\mu_K\ar[r]\ar[d]& 0 \\ 
0\ar[r]& G^\mu_K\ar[r]& G^{\rig}(\rho)[p^\infty]\ar[r]\ar[d]& G^{\rig}(\rho)[p^\infty]/G^\mu_K\ar[r]\ar[d]& 0 \\ 
& & G^\et_K\ar@{=}[r]\ar[d]& G^\et_K\ar[d] \\ & &0&0}
\end{equation}
with exact rows and columns.

\textbf{Claim:} $G^{\rig}(\rho)[p^\infty]/G^\mu_K$ extends to a unique $p$-divisible group over $\cO_K$ up to unique isomorphism. 

The uniqueness follows from Tate's theorem (see \cite[Thm. 4]{tat1}). We show the existence. Let $H_K$ be the Cartier dual of $G^{\rig}(\rho)[p^\infty]/G^\mu_K$. It suffices to show that $H_K$ extends a $p$-divisible group over $\cO_K$. Let $\tau$ denote the Galois $\Z_p$-representation associated to $H_K$. Then $\tau$ as a subrepresentation of $\rho^\vee(1)$ is semistable with Hodge-Tate weights in $\{0,1\}$ by \cite[Cor. 3.16 (ii)]{bri1}, where $(-)^\vee$ denotes the dual representation and $(-)(1)$ denotes the Tate twist of Galois representation. Since $G^\circ/G^\mu$ is connected and has no multiplicative part, its Cartier dual is connected and thus the associated Galois representation has no non-trivial potentially unramified quotient. Obviously the Galois representation associated to the Cartier dual $(G^\et)^*$ of $G^\et$ has no non-trivial potentially unramified quotient. It follows that $\tau$ has no non-trivial potentially unramified quotient. Then the $p$-divisible rigid analytic group associated to $\tau$ is connected and comes from a $p$-divisible formal group $F_{\tau}$ over $\cO_K$ by \cite[Prop. 2.1]{Far22} and \cite[Thm. 6.1]{Far19}. Thus we have $H_K=(F_{\tau}[p^\infty])_K$, i.e. the $p$-divisible group $H:=F_{\tau}[p^\infty]$ over $\cO_K$ extends $H_K$. This finishes the proof of the claim.

 Now the upper row and the right-most column of (\ref{pan-ext diagram associated a sst represenation of HT weights 0 and 1}) together extend  to the following diagram
\[\xymatrix@R=10pt{& & & 0\ar[d] & \\ 
0\ar[r] & G^\mu\ar[r]& G^\circ\ar[r]^-{\alpha}\ar[rd]_{\gamma} & G^\circ/G^\mu\ar[r]\ar[d]^{\beta} & 0 \\ 
& & & H^*\ar[d] \\ 
& & & G^\et\ar[d] \\ 
& & &0}\]
with exact row and column over $\cO_K$ and $\gamma=\beta\circ\alpha$. To finish the proof, we need to show that $\ker(\gamma)$ (resp. $\coker(\gamma)$) is a multiplicative (resp. \'etale) $p$-divisible group. But this is clear, as $\ker(\gamma)=\ker(\alpha)$ and $\coker(\gamma)=\coker(\beta)$.
\end{proof}

\begin{cor}\label{sst rep. with HT weights 0 and 1 arise from log BT}
 Any $\rho\in \Repz^{\sst, \{0,1\}}(\Galk)$ arises from a dual representable log $p$-divisible group under the functor $T_p\colon \btlogd\to \Repz(\Galk)$.
\end{cor}
\begin{proof}
 This follows from Proposition \ref{the $p$-div gp associated to a sst rep. with HT weights 0 and 1 is sst} and Theorem \ref{BZ}.
\end{proof}

\subsection{Proof of the second part of Theorem \ref{main}}\label{s.Bb}
%%%%%%%%%%%%%%%%%%%%%%%%%

Now we are ready to prove the second part of Theorem \ref{main}. We make it into a separate theorem.

\begin{thm}\label{BZb}
 The functor 
 \[T_p\colon \btlogd\to {\bf Rep}_{\Z_p}(\Galk) \text{ (resp. $V_p\colon \btlogd\otimes\Q\to \Repq(\Galk)$)} \]
 is fully faithful and has essential image the full subcategory $\Repz^{\sst,\{0,1\}}(\Galk)$ (resp. $\Repq^{\sst,\{0,1\}}(\Galk)$).
 %in $ {\bf Rep}_{\Z_p}(\Galk)$.
\end{thm}
\begin{proof}
We only prove the case of $T_p$. By definition, the functor $T_p$ in the statement is the composition $\btlogd\xrightarrow{(\ )_K} {\bf BT}_{K}\xrightarrow{T_p}{\bf Rep}_{\Z_p}(\Galk)$. By Theorem \ref{BZ}, the functor $(\ )_K$ is fully faithful. Since the base field $K$ is of characteristic 0, the functor $T_p\colon {\bf BT}_{K}\to {\bf Rep}_{\Z_p}(\Galk)$ is fully faithful. It follows that the functor $T_p$ in the statement is fully faithful.  The rest follows from Proposition \ref{the rep. associated to a log BT is sst with HT weight 0 and 1} and Corollary \ref{sst rep. with HT weights 0 and 1 arise from log BT}.
\end{proof}

\subsection{Consequences of Theorem \ref{main}}\label{section.corollaries}
%%%%%%%%%%%%%%%%%%%%%%%%%

The \emph{proof of Corollary \ref{c.NOS}} is now immediate. 
 By Theorem \ref{main}, it suffices to show that $A_K$ has semistable reduction if and only if $A_K[p^\infty]$ has semistable reduction,  and this follows from \cite[Expos\'e IX]{sga72} as explained in \cite[proof of 2.5]{dJong98}.

Note that this corollary was already proved in \cite[Cor. 5.3.4]{Bre00}, with conditions on $K$, and the general case can be deduced from results in \cite{Kis06} similarly.

 Also the \emph{proof of Corollary \ref{c.stronglydivisible}} is now easy. 
By \cite[Conj. 2.2.6 (1)]{Bre02} (proved in  \cite[Thm. 2.3.5]{Liu08} for $p>2$) and \cite[Thm. 2.2.7 (2)]{Bre02}, there is an anti-equivalence between the category of strongly divisible modules of weight $\leq1$ and the category ${\bf Rep}^{\mr{st},\{0,1\}}_{\Z_p}(\Galk)$. Then one concludes by  Theorem \ref{main}(b).

\subsection{Compatibility of Kato monodromy and Fontaine monodromy}\label{s.compatibilityKF}
%%%%%%%%%%%%%%%%%%%%%%%%%

The diagram \eqref{triangle_diagram_main3} extends to a diagram
\[\xymatrix@C=20pt{
&\btlogd  \ar[ld]_-{(\ )_K}^-{\simeq}\ar[rd]^-{T_p}_-{\simeq}  \\
{\bf BT}_{K}^{\mr{st}}\ar[rr]^-{\simeq}_-{T_p}  &&{\bf Rep}^{\mr{st},\{0,1\}}_{\Z_p}(\Galk)\ar[r]_-{\otimes_{\Z_p}\Q_p} &{\bf Rep}^{\mr{st},\{0,1\}}_{\Q_p}(\Galk)\ar[r]^-{\simeq}_-{D_{\sst}} &\underline{M}^{\mathrm{a},\{-1,0\}}
},\]
where 
\begin{enumerate}[(1)]
\item ${\bf Rep}^{\mr{st},\{0,1\}}_{\Q_p}(\Galk)$ denotes the category of semistable $\Q_p$-representation of $\Galk$ with Hodge-Tate weights in $\{0,1\}$,
\item $\underline{M}^{\mathrm{a}}$ denotes the category of admissible filtered $(\varphi,N)$-modules over $K$ (see \cite[\S 4.1]{CF00}), and $\underline{M}^{\mathrm{a},\{-1,0\}}$ denotes the full subcategory of $\underline{M}^{\mathrm{a}}$ consisting of objects $D$ such that $\mathrm{Fil}^{-1}D_K=D_K$ and $\mathrm{Fil}^{1}D_K=0$,
\item the functor $D_{\sst}\colon \Repq^{\sst}(\Galk)\xrightarrow{\simeq}\underline{M}^{\mathrm{a}}, V\mapsto (V\otimes_{\Q_p}B_{\sst})^{\Galk}$  (see \cite[\S 4.1]{CF00}, as well as Subsection \ref{Subsection.pdiv-to-crysrep}) is an exact tensor functor, as well as an equivalence of categories by \cite[Prop. 4.2]{CF00}. Here we take its restriction to $\Repq^{\sst,\{0,1\}}(\Galk)$, which is an equivalence with the category $\underline{M}^{\mathrm{a},\{-1,0\}}$. 
\end{enumerate}

For any object in $\btlogd$ (resp. in ${\bf BT}_{K}^{\mr{st}}$), there is an associated Kato monodromy (resp. Grothendieck monodromy). By Theorem \ref{thm.monodromypairings} Kato monodromy is compatible with Grothendieck monodromy along the equivalence of categories $(\ )_K$. It seems that there is no monodromy associated to a representation $\rho$ in ${\bf Rep}^{\mr{st},\{0,1\}}_{\Z_p}(\Galk)$. Nevertheless if we pass to $\underline{M}^{\mathrm{a},\{-1,0\}}$, we have the $K_0$-linear endomorphism $N$ on the $K_0$-vector space $D_{\sst}(\rho\otimes_{\Z_p}\Q_p)$. 
Given an object $H\in\btlogd$, it is natural to investigate the relation between the Kato monodromy $\beta$ of $H$ and the map $N$ on $D_{\sst}(V_p(H))$. This subsection is devoted to this investigation. 

For convenience, for $V\in {\bf Rep}^{\mr{st},\{0,1\}}_{\Q_p}(\Galk)$ we call the map $N$ of $D_{\sst}(V)$ the \emph{Fontaine monodromy map} of $V$, as well as of $D_{\sst}(V)$ by abuse of terminology.

 Let $0\to H^\circ\to H\to H^\et\to 0$ be the connected-\'etale decomposition (see \eqref{e.connected-etale}) of $H$, and let $(H^{\cl},\beta)$ be the pair guaranteed by Theorem \ref{thm1.2}, and thus $\beta\colon H^\et(1)\to H^\circ$ is the Kato monodromy of $H$ according to Definition \ref{def.katomonodromy}. 
 Let $V_\pi\in {\bf Rep}^{\mr{st},\{0,1\}}_{\Q_p}(\Galk)$ be the representation \eqref{the rep. associated to the log 1-mot M_pi} associated to the log 1-motive $M_\pi=[\Z\stackrel{1\mapsto \pi}{\longrightarrow} \Gml]$. It fits into  a short exact sequence $0\to \Q_p(1)\xrightarrow{a} V_{\pi}\xrightarrow{b} \Q_p\to 0$ in ${\bf Rep}^{\mr{st},\{0,1\}}_{\Q_p}(\Galk)$, see \eqref{eq1}.
 Since the Fontaine monodromy map decreases the slope by $1$ (see \cite[\S3.3]{CF00}), the Fontaine monodromy map of $V_\pi$ factors as
\[D_{\sst}(V_\pi)\xrightarrow{D_{\sst}(b)} D_{\sst}(\Q_p)\xrightarrow{\overline{N}_{\pi}} D_{\sst}(\Q_p(1))\xrightarrow{D_{\sst}(a)} D_{\sst}(V_{\pi}).\]
Since $D_{\sst}\colon {\bf Rep}^{\mr{st}}_{\Q_p}(\Galk)\xrightarrow{\simeq}\underline{M}^{\mathrm{a}}$ is a tensor functor, the Fontaine monodromy map of $V_{p}(M_{\pi}[p^\infty]\otimes_{\Z_p}H^\et)$ factors as  
\begin{align*}
D_{\sst}(V_{p}(M_{\pi}[p^\infty]\otimes_{\Z_p}H^\et))\to D_{\sst}(V_{p}(H^\et))\xrightarrow{1_{D_{\sst}(V_{p}(H^\et))}\otimes\overline{N}_{\pi}} D_{\sst}(V_{p}(H^\et)(1))  \\
\to D_{\sst}(V_{p}(M_{\pi}[p^\infty]\otimes_{\Z_p}H^\et)).
\end{align*}

 Let $N_{\beta}$ be the composition
\begin{multline}\label{e.Nbeta}
N_{\beta}\colon D_{\sst}(V_{p}(H))\to D_{\sst}(V_{p}(H^\et))\xrightarrow{1_{D_{\sst}(V_{p}(H^\et))}\otimes\overline{N}_{\pi}}  D_{\sst}(V_{p}(H^\et)(1)) \\
\xrightarrow{D_{\sst}(V_{p}(\beta))} D_{\sst}(V_p(H^\circ))\to  D_{\sst}(V_{p}(H)).
\end{multline}
Note that the operator $N_\beta$ has mixed information from Fontaine's monodromy and Kato's monodromy.

The main result of this subsection is the following theorem.

\begin{thm}\label{thm.monodromyKatoFontaine}
 The Fontaine monodromy map of $D_{\sst}(V(H))$ agrees with $N_{\beta}$  defined in \eqref{e.Nbeta}.
\end{thm}

 We need a lemma for proving it.

\begin{lem}\label{reduced monodromy map}
 Let $D_1,D_2\in \underline{M}^{\mathrm{a},\{-1,0\}}$. 
\begin{enumerate}[(a)]
\item The slopes of $D_1$ and $D_2$ lie in $[-1,0]$.
\item Assume that $D_1$ is of slope $0$ and the slopes of $D_2$ lie in $[-1,0)$. Then 
\begin{enumerate}[(b.1)]
    \item Both $D_1$ and $D_2$ have trivial Fontaine monodromy map.
    \item Let $0\to D_2\xrightarrow{i} D\xrightarrow{q} D_1\to 0$ be an extension in the abelian category $\underline{M}^{\mathrm{a},\{-1,0\}}$. Then the Fontaine monodromy map $N_D$ of $D$ factors as
    \[D\xrightarrow{q} D_1\xrightarrow{\overline{N}_{D}} D_2\xrightarrow{i} D.\]
    \item Let $f:D_2\to D_3$ be a map in $\underline{M}^{\mathrm{a},\{-1,0\}}$. Assume that the slopes of $D_3$ lie in $[-1,0)$. Let $f_*D$ be the pushout of $D$ along $f$, as depicted in the diagram
    \[\xymatrix{
    0\ar[r] &D_2\ar[r]^i\ar[d]^f &D\ar[r]^q\ar[d]^a &D_1\ar[r]\ar@{=}[d] &0 \\
    0\ar[r] &D_3\ar[r]^j &f_*D\ar[r]^c &D_1\ar[r] &0
    }.\]
    Let $\overline{N}_{f_*D}$ be the map defined in  (b.2) for $f_*D$. Then 
    \[\overline{N}_{f_*D}=f\circ\overline{N}_D.\]
    \item Let $g\colon D_4\to D_1$ be a map in $\underline{M}^{\mathrm{a},\{-1,0\}}$. Assume that $D_4$ is of slope 0. Let $g^*D$ be the pullback of $D$ along $g$, as depicted in the diagram 
    \[\xymatrix{
    0\ar[r] &D_2\ar[r]^i &D\ar[r]^q &D_1\ar[r] &0 \\
    0\ar[r] &D_2\ar@{=}[u]\ar[r]^d &g^*D\ar[r]^e\ar[u]^b &D_4\ar[r]\ar[u]^{g} &0
    }.\]
    Let $\overline{N}_{g^*D}$ be the map defined in  (b.2) for $g^*D$. Then 
    \[\overline{N}_{g^*D}=\overline{N}_D\circ g.\]
    \item Let $0\to D_2\xrightarrow{i'} D'\xrightarrow{q'} D_1\to 0$ be another extension of $D_1$ by $D_2$, and let $D+_{\mathrm{B}}D'$ be the Baer sum of extensions in $\underline{M}^{\mathrm{a},\{-1,0\}}$. Then we have
    \[\overline{N}_{D+_{\mathrm{B}}D'}=\overline{N}_D+\overline{N}_{D'}.\]
    \end{enumerate}
\end{enumerate}
\end{lem}

\begin{proof}
(a) This follows from the admissibility.

 We prove the assertions in (b). Both (b.1) and (b.2) follow from the fact that the Fontaine monodromy map decreases the slope by 1.

 (b.3) Since $f_*D$ is a quotient of $D_3\oplus D$, the Fontaine monodromy map $N_{f_*D}$ of $f_*D$ is induced by that of $D_3\oplus D$. We have the following commutative diagram
\[\xymatrix{
&D_3\oplus D\ar[d]^(.25){(0,N_D)}\ar[rr]^-{(0,q)}\ar[ld]_{(j,a)} &&D_1\ar[d]^{\overline{N}_D} \\
f_*D\ar[d]_{N_{f_*D}}\ar[rrru]|(.37)\hole ^(.25)c &D_3\oplus D\ar[ld]^{(j,a)} &&D_2\ar[ll]_-{(0,i)}\ar[ld]^{f} \\
f_*D & &D_3\ar[ll]^j
},\]
from which we get $j\circ (f\circ\overline{N}_D)\circ c=N_{f_*D}$.
Since $c$ is surjective and $j$ is injective, we must have $\overline{N}_{f_*D}=f\circ\overline{N}_D$.

(b.4) Since $g^*D$ is a subobject of $D\oplus D_4$, the Fontaine monodromy map $N_{g^*D}$ of $g^*D$ is induced by that of $D\oplus D_4$ which is $N_D\oplus0$. We have the following commutative diagram
\[\xymatrix{
g^*D\ar[d]_{N_{g^*D}}\ar[rd]^{(b,e)}\ar[rr]^e &&D_4\ar[rd]^{g} \\
g^*D\ar[rd]_{(b,e)} &D\oplus D_4\ar[rr]^{(q,0)}\ar[d]^(.25){N_D\oplus0} &&D_1\ar[d]^{\overline{N}_D} \\
 &D\oplus D_4 &&D_2\ar[ll]^{(i,0)}\ar[lllu]|(.64)\hole_(.3)d \\
}\]
from which we get $N_{g^*D}=d\circ(\overline{N}_D\circ g)\circ e$. Since $e$ is surjective and $d$ is injective, we must have $\overline{N}_{g^*D}=\overline{N}_D\circ g$.

(b.5) By definition $D+_{\mathrm{B}}D'$ is constructed as in the following diagram
\[\xymatrix{
0\ar[r] &D_2\oplus D_2\ar[r]^{i\oplus i'}\ar[d]_{+} &D\oplus D'\ar[r]^{q\oplus q'}\ar[d] &D_1\oplus D_1\ar[r]\ar@{=}[d] &0\\
0\ar[r] &D_2\ar[r] &E\ar[r] &D_1\oplus D_1\ar[r] &0\\
0\ar[r] &D_2\ar[r]\ar@{=}[u] &D+_{\mathrm{B}}D'\ar[r]\ar[u] &D_1\ar[r]\ar[u]_{\Delta} &0
},\]
where $E:=+_*(D\oplus D')$ denotes the pushout of $D\oplus D'$ along the sum map $+\colon D_2\oplus D_2\to D_2$ and $D+_{\mathrm{B}}D'$ is the pullback of $E$ along the diagonal map $\Delta\colon D_1\to D_1\oplus D_1$. By (b.3) and (b.4), we get 
\[\overline{N}_{D+_{\mathrm{B}}D'}=+\circ(\overline{N}_D\oplus\overline{N}_{D'})\circ\Delta=\overline{N}_D+\overline{N}_{D'}.\]
\end{proof}

\begin{proof}[Proof of Theorem \ref{thm.monodromyKatoFontaine}]
For any $G\in\btlogd$, let $N_{G}$ be the Fontaine monodromy map of $D_G:=D_{\sst}(V_p(G))=D_{\sst}(T_p(G)\otimes_{\Z_p}\Q_p)$.

Now let $H\in\btlogd$ be as in the beginning of this subsection. Let $H^\beta:=\Phi_2(\beta)$, and thus $H$ is the Baer sum of $H^\cl$ and $H^\beta$ as extensions of $H^\et$ by $H^\circ$. 
We consider $D_H$ (resp. $D_{H^\cl}$, resp. $D_{H^\beta}$) as an extension of $D_{H^\et}$ by $D_{H^\circ}$ canonically, and let $\overline{N}_{H}$ (resp. $\overline{N}_{H^\cl}$, resp. $\overline{N}_{H^\beta}$) be the map $D_{H^\et}\to D_{H^\circ}$ defined as in Lemma \ref{reduced monodromy map} (b.2). By construction, $D_H$ is the Baer sum of $D_{H^\cl}$ and $D_{H^\beta}$ as extensions of $D_{H^\et}$ by $D_{H^\circ}$. By Lemma \ref{reduced monodromy map} (b.5), we have $\overline{N}_{H}=\overline{N}_{H^\cl}+\overline{N}_{H^\beta}=\overline{N}_{H^\beta}$ which can be further computed as $D_{\sst}(V_p(\beta))\circ \overline{N}_{H^\et\otimes_{\Z_p}M_{\pi}[p^\infty]}$ by Lemma \ref{reduced monodromy map} (b.3). Previously we have seen that $\overline{N}_{H^\et\otimes_{\Z_p}M_{\pi}[p^\infty]}=1_{D_{\sst}(V_{p}(H^\et))}\otimes\overline{N}_{\pi}$. Therefore, $N_H$ agrees with $N_\beta$.  
\end{proof}

Theorem \ref{thm.monodromyKatoFontaine} tells us that the Fontaine monodromy map $N_H$ of $D_{\sst}(V_{p}(H))$ is determined by the Kato monodromy map of $H$, i.e. one direction of Theorem \ref{t.KFmonodromy} holds. Now we discuss the other direction. Using \cite[\S II.4]{berg1}, one can easily compute $\overline{N}_{\pi}$ which is an isomorphism of $K_0$-vector spaces. It follows that 
\[V_{p}(\beta)=V_{\sst}\left( \overline{N}_{H}\circ (1_{D_{\sst}(V(H^\et))}\otimes\overline{N}_{\pi})^{-1}\right),\]
where $V_{\sst}\colon \underline{M}^{\mathrm{a}}\xrightarrow{\simeq}{\bf Rep}^{\sst}_{\Q_p}(\Galk)$ is the functor defined in the third paragraph of \cite[\S4.1]{CF00} which is quasi-inverse to $D_{\sst}$ by \cite[\S4.1 Cor.]{CF00}. 
One should not confuse $V_{\sst}$ with $V_p$ which associates to an object of $\bt_{\cO_K}$ (or $\btlogd$) its $\Q_p$-representation of $\Galk$. Since $V_p\colon \bt_{\cO_K}\otimes\Q\to \Repq^{\cris,\{0,1\}}(\Galk)$ is an equivalence of categories by Theorem \ref{classicalA}, the Kato monodromy $\beta$ of $H$ is rationally determined by $N_H$. 

%%%%%%%%%%%%%%%%%%%%%%%%%
\section*{Acknowledgements}
%%%%%%%%%%%%%%%%%%%%%%%%%

The authors wish to thank an anonymous referee for bringing their attention to Corollary \ref{c.stronglydivisible} as well as Professor Christophe Breuil and Professor Tong Liu for answering a question on it.
The first author was partially supported by Italian  PRIN 2022 ``The arithmetic of motives and L-functions'', project number 20222B24AY, and by the project BIRD228095   "Arithmetic cohomology theories" of University of Padova. 
Part of this article was written during S. Wang's visit to Westlake University in January 2023. He would like to thank Professor Yigeng Zhao for his invitation and hospitality. Shanwen Wang is supported by the Fundamental Research Funds for the Central Universities, and the Research Funds of Renmin University of China No.20XNLG04 and The National Natural Science Foundation of China (Grant No.11971035).
Part of this article was written when the third author was a member of the research group of Professor Ulrich G\"ortz under partial support from Research Training Group 2553 of the German Research Foundation DFG. He wishes to thank Professor Xu Shen for a very helpful discussion on the $p$-adic criterion for semistable reduction of abelian varieties.

\bibliographystyle{siam}

\end{document}